\documentclass [leqno] {amsart}
  
\usepackage {amssymb}
\usepackage {amsthm}
\usepackage {amscd}

\usepackage[all]{xy}
\usepackage {hyperref} 

\begin{document}
\bibliographystyle{alpha}
\theoremstyle{plain}
\newtheorem{proposition}[subsection]{Proposition}
\newtheorem{lemma}[subsection]{Lemma}
\newtheorem{corollary}[subsection]{Corollary}
\newtheorem{thm}[subsection]{Theorem}
\newtheorem{introthm}{Theorem}
\newtheorem*{thm*}{Theorem}
\newtheorem{conjecture}[subsection]{Conjecture}
\newtheorem{question}[subsection]{Question}
\newtheorem{fails}[subsection]{Fails}

\theoremstyle{definition}
\newtheorem{definition}[subsection]{Definition}
\newtheorem{notation}[subsection]{Notation}
\newtheorem{condition}[subsection]{Condition}
\newtheorem{example}[subsection]{Example}
\newtheorem{claim}[subsection]{Claim}

\theoremstyle{remark}
\newtheorem{remark}[subsubsection]{Remark}

\numberwithin{equation}{subsection}

\newcommand{\eq}[2]{\begin{equation}\label{#1}#2 \end{equation}}
\newcommand{\ml}[2]{\begin{multline}\label{#1}#2 \end{multline}}
\newcommand{\mlnl}[1]{\begin{multline*}#1 \end{multline*}}
\newcommand{\ga}[2]{\begin{gather}\label{#1}#2 \end{gather}}
\newcommand{\mat}[1]{\left(\begin{smallmatrix}#1\end{smallmatrix}\right)}

\newcommand{\arir}{\ar@{^{(}->}}
\newcommand{\aril}{\ar@{_{(}->}}
\newcommand{\are}{\ar@{>>}}

\newcommand{\xr}[1] {\xrightarrow{#1}}
\newcommand{\xl}[1] {\xleftarrow{#1}}
\newcommand{\lra}{\longrightarrow}
\newcommand{\inj}{\hookrightarrow}

\newcommand{\mf}[1]{\mathfrak{#1}}
\newcommand{\mc}[1]{\mathcal{#1}}
\newcommand{\mbf}[1]{\mathbf{#1}}

\newcommand{\CH}{{\rm CH}}
\newcommand{\Gr}{{\rm Gr}}
\newcommand{\codim}{{\rm codim}}
\newcommand{\cd}{{\rm cd}}
\newcommand{\Spec} {{\rm Spec}}
\newcommand{\supp} {{\rm supp}}
\newcommand{\Hom} {{\rm Hom}}
\newcommand{\End} {{\rm End}}
\newcommand{\id}{{\rm id}}
\newcommand{\Aut}{{\rm Aut}}
\newcommand{\sHom}{{\rm \mathcal{H}om}}
\newcommand{\Tr}{{\rm Tr}}
\newcommand{\Ext}{{\rm Ext}}


\renewcommand{\P} {\mathbb{P}}
\newcommand{\Z} {\mathbb{Z}}
\newcommand{\Q} {\mathbb{Q}}
\newcommand{\C} {\mathbb{C}}
\newcommand{\F} {\mathbb{F}}

\newcommand{\OO}{\mathcal{O}}

 
\title{On integrality of $p$-adic iterated integrals}

\author{Andre Chatzistamatiou}
\email{a.chatzistamatiou@uni-due.de}

\thanks{The author is supported by the Heisenberg program (DFG)}

\begin{abstract}
The purpose of this paper is to prove integrality for certain $p$-adic iterated Coleman integrals. As underlying geometry we will take the complement of a divisor $D\subset X$ with good reduction, where $X$ is the projective line or an elliptic curve over the Witt vectors of a perfect characteristic $p$ field. As a corollary we prove a lower bound for the valuations of $p$-adic multiple zeta values. 
\end{abstract}

\maketitle


\section*{Introduction}
The basic idea for analytic continuation in real analytic spaces fails for classical $p$-adic spaces due to their totally disconnected nature. On such spaces it is not at all evident how to patch local primitives of a function together into an essentially unique global primitive. Guided by Dwork's principle, that the $p$-adic analogue of analytic continuation along a path is "analytic continuation along Frobenius", Coleman developed in  \cite{C1,C2,CdS} a theory for solving iterated integrals. The solutions form a subring of the locally analytic functions whose elements are called Coleman functions. As an application, Coleman constructed the $p$-adic version of the polylogarithm functions. The relationship between values of "the" $p$-adic zeta function at positive integers and the value of Coleman polylogarithms at $1$ is the same as in the complex setting.       

In our paper, we will use the generalization of Coleman's theory by Besser \cite{B} and, independently, Vologodsky \cite{V}, which is built on the Tannakian formalism. The notion of a set of paths between two points is replaced by the space of natural isomorphisms between two fibre functors on a suitable category of unipotent connections. The Frobenius action on the space of natural isomorphisms singles out a distinguished path, which is called the \emph{Frobenius invariant path}. Therefore "analytic continuation along Frobenius" admits a precise formulation as continuation along the Frobenius invariant path. 

The main tool in our paper is an extension of the Tannakian formalism to very simple unipotent categories over any commutative ring (Corollary \ref{corollary-Tanaka}). Broadly speaking, these categories are equivalent to unipotent representations of a free tensor algebra with finitely many generators. We develop the theory in the first section to a degree sufficient for our application; a more general study of Tannakian duality in the unipotent case would be very worthwhile. Our application is to the category $\mc{U}$ of unipotent connections on $\hat{X}$ with logarithmic singularities in $D\neq \emptyset$, where $X$ is a geometrically connected smooth projective curve over the ring of Witt vectors $W$ of a perfect characteristic $p$ field, $\hat{X}$ denotes the completion along the special fibre, and $D$ is finite \'etale over $W$.    

In Sections \ref{section-unipotent-connections} and \ref{section-Frobenius-fibre-functors}, we recall some properties of $\mc{U}$. In particular, we describe the absolute Frobenius functor, and we give integral versions (that is, defined over $W$) of the fibre functors attached to $W$-points of $X$ with good or bad reduction with respect to $D$, and give an account of tangential fibre functors.  

Section \ref{section-main-thms} contains the two main theorems. Theorem \ref{thm-reduction-to-projective-curve} asserts that it is sufficient to construct an integral version of the Frobenius invariant path between two fibre functors $x$ and $y$ for unipotent connections on $\hat{X}$. Then there is a unique extension to connections with logarithmic singularities. In the case $X=\P^1_W$, this implies integrality immediately. For an elliptic curve $X$, we can only prove integrality under the additional assumption that the reductions $x_0,y_0\in X(k)$ of the points underlying $x,y$, yield a point $x_0-y_0$ of order prime to $p$ (Theorem \ref{thm-construction-path}). 

In Section \ref{section-application-mutliple-zeta-values}, we give the application to $p$-adic multiple zeta values (Corollary \ref{corollary-multiple-zeta-values}): $\zeta_p(k_1,\dots,k_m)\in p^{\left[ \sum_{i=1}^m k_i\right]}$, for all positive integers $k_1,\dots,k_m$, all primes $p$, and where $p^{\left[ \sum_{i=1}^m k_i\right]}\subset \Z_p$ is the ideal generated by all $\frac{p^j}{j!}$ with $j\geq \sum_{i=1}^m k_i$. 

\section{A simple case of unipotent Tannaka duality for arbitrary base rings}

Let $R$ be a commutative ring, and let $\mc{C}$ be an $R$-linear exact category. Let $\mbf{1}$ be an object of $\mc{C}$ such that $\Hom(\mathbf{1},\mathbf{1})=R$. 

\begin{definition}\label{definition-Ur}
	For any non-negative integer $r$, we define $\mc{U}_r=\mc{U}_r(\mc{C},\mbf{1})$ to be the full subcategory consisting of objects $V\in \mc{C}$ that have  an admissible filtration (that is, the inclusions are admissible morphisms)
	$$
	0\subset V_0\subset V_1 \subset \dots \subset V_r=V,
	$$
	such that the quotients $V_{i+1}/V_i$ and $V_0$ are isomorphic to a finite sum with summands $\mathbf{1}$. We set $\mc{U}=\bigcup_r \mc{U}_r$. 
\end{definition}

Note that $\mc{U}$ is closed under extensions in $\mc{C}$ and is an $R$-linear exact category. Let $\mc{U}_c$ be the idempotent completion of $\mc{U}$. We have its derived category $D(\mc{U}_c)$ at disposal \cite{Neeman1990}, and may set 
$$
\Ext^i(X,Y):=\Hom_{D(\mc{U}_c)}(X,Y[i]),
$$
for $X,Y\in \mc{U}$. As for abelian categories, this construction is the same as the Yoneda construction where an element in $\Ext^i(X,Y)$ is represented by an exact sequence of admissible morphisms in $\mc{U}_c$, which starts with $Y$ and ends with $X$.  

\begin{notation}\label{notation-free}
	We say that $\mc{U}$ is \emph{free} if $\Ext^1(\bf{1},\bf{1})$ is a free and finitely generated $R$-module and ${\rm Ext}^2(\mathbf{1},\mbf{1})=0$.
\end{notation}

\begin{definition}\label{definition-fibre-functor}
	Let $r\geq 0$. An $R$-linear functor $\omega:\mc{U}_r\xr{} \text{($R$-modules)}$ is called a \emph{fibre functor} if $\omega$ respects exact sequences and $\omega(\mbf{1})\cong R$. We define fibre functors with source category $\mc{U}$ in the same way.
\end{definition}

By definition, we have $\mc{U}_0\cong \text{(free finitely generated $R$-modules)}$. For a free finitely generated $R$-module $V$, we will denote by $V\otimes \mbf{1}\in \mc{U}_0$ the representing object of the functor $\mc{U}_0 \xr{} \text{(free finitely generated $R$-modules)}; T\mapsto \Hom(T,\mbf{1})\otimes_R V$. It comes equipped with an isomorphism $\Hom(\mbf{1},V\otimes \mbf{1})\xr{} V$.

\begin{proposition}\label{proposition-representation}
	Suppose that $\mc{U}$ is free, and let $r\geq 0$. 
	\begin{enumerate}	
		\item There is an object $E_r$ of $\mc{U}_r$ such that $\Hom(E_r,-):\mc{U}_r\xr{} \text{($R$-modules)}$ is a fibre functor. 
		\item Let $E_r$ be an object with the properties of (1), let $\epsilon:E_r\xr{} \mbf{1}$ be a generator for $\Hom(E_r,\mbf{1})$. For every fibre functor 
		$\omega:\mc{U}_r\xr{} \text{($R$-modules)}$, and every $e\in \omega(E_r)$ such that $\omega(\epsilon)(e)$ generates $\omega(\mbf{1})$, $\omega$ is represented by $(E_r,e)$.	
	\end{enumerate}
	\begin{proof}
		For $r=0$, we can take $E_0=\mbf{1}$. 
		
		Set $\Gamma:= {\rm Ext}^1(\mathbf{1},\mbf{1})^{\vee}$ and $T^{\otimes r}:= \Gamma^{\otimes r}\otimes \mbf{1}$. For $r>0$, we will proceed by induction and construct an extension 
		\begin{equation}\label{equation-representation-extension}
		0\xr{} T^{\otimes r} \xr{} E_r \xr{} E_{r-1} \xr{} 0,
		\end{equation}
	such that the induced map 
	\begin{equation}\label{equation-representation-induction}
		\Hom(T^{\otimes r},\mbf{1})\xr{} \Ext^1(E_{r-1},\mbf{1})	
	\end{equation}
	is an isomorphism. Note that in view of the embedding theorem \cite[A.7.1,A.7.16]{TT} (cf.\cite[Theorem~A.1]{Buehler}), the category $\mc{U}$ is closed under extensions in $\mc{U}_c$. 
        
	Obviously, ${\rm Ext}^2(\mathbf{1},\mbf{1})=0$ implies ${\rm Ext}^2(V,W)=0$ for any two objects $V,W\in \mc{U}_r$. Suppose we have constructed $E_{r-1}$. Then 
	$$
	\Ext^1(E_{r-1},\mbf{1}) \xr{\cong} \Ext^1(T^{\otimes {r-1}},\mbf{1})
	$$
	is an isomorphism.
	Let $E_r$ correspond to the preimage  of the identity via the map
	\begin{multline*}
	\Ext^1(E_{r-1},T^{\otimes r}) \xr{} \Ext^1(T^{\otimes {r-1}},T^{\otimes r}) \xr{}  \Ext^1(\mbf{1},\mbf{1})\otimes \Gamma^{\vee \otimes r-1}  \otimes \Gamma^{\otimes r} \\ =\Hom(\Gamma^{\otimes r},\Gamma^{\otimes r}). 
	\end{multline*}
	Let us show \eqref{equation-representation-induction} holds. This map is induced by the class of $E_r$ and the first row of the commutative diagram
	$$
	\xymatrix{
		\Hom(T^{\otimes r},\mbf{1})\otimes  \Ext^1(E_{r-1},T^{\otimes r}) \ar[r] \ar[d]^{\cong}& \Ext^1(E_{r-1},\mbf{1}) \ar[d]^{\cong} \\
		\Hom(T^{\otimes r},\mbf{1})\otimes  \Ext^1(T^{\otimes {r-1}},T^{\otimes r}) \ar[r] & \Ext^1(T^{\otimes {r-1}},\mbf{1}),
	}
	$$
	which proves the claim. 
	
	By induction, we easily obtain $\Hom(E_r,\mbf{1})\cong R$. 
	
	Let us prove by induction on $r$ that 
	\begin{equation} \label{equation-representation-Er-easy-sequence}
	0\xr{} \Hom(E_r,V_{0}) \xr{} \Hom(E_r,V) \xr{} \Hom(E_r,V_{r-1})\xr{} 0
	\end{equation}  
	is an exact sequence if $V\in \mc{U}_r$, and $V_0\subset V$ is such that $V_0\in \mc{U}_0$ and $V_{r-1}:=V/V_0\in \mc{U}_{r-1}$.
	
	First, we note that, for every $V\in \mc{U}_{r-1}$, the map  
	$$
	\Hom(E_{r-1},V)\xr{} \Hom(E_{r},V),
	$$
	induced by $E_r\xr{} E_{r-1}$, is an isomorphism. Certainly, this holds for $V\in \mc{U}_0$. For the general case, we can use an exact sequence
	$$
	0\xr{} V_{0} \xr{} V \xr{} V_{r-2} \xr{} 0,
	$$
	with $V_0\in \mc{U}_0, V_{r-2}\in \mc{U}_{r-2}$, and induction.

        For  $\eta\in \Hom(E_{r-1},V_{r-1})=\Hom(E_r,V_{r-1})$ 
   	consider 
	$$
	\xymatrix{
		0 \ar[r] & V_0 \ar[r] & V \ar[r] & V_{r-1} \ar[r] & 0 \\
		0 \ar[r] & V_0 \ar[r]\ar[u]^{=} & E_{\eta} \ar[r] \ar[u] & E_{r-1} \ar[r] \ar[u]^{\eta} &  0, 
	}
	$$
	with $E_{\eta}$ being the induced extension. Since $\Hom(T^{\otimes r},V_0)\xr{} \Ext^1(E_{r-1},V_0)$ is an isomorphism, $E_{\eta}$ defines a morphism $\tau:T^{\otimes r}\xr{} V_0$, and the class 
	of $E_{\eta}$ is the image of the class of $E_{r}$ via $\tau$. We obtain a commutative diagram
       	$$
		\xymatrix{
			0 \ar[r] & V_0 \ar[r] & E_{\eta} \ar[r]  & E_{r-1} \ar[r]  &  0,\\
			0 \ar[r] & V_0 \ar[r]\ar[u]^{\tau} & E_{r} \ar[r] \ar[u]^{\psi} & E_{r-1} \ar[r] \ar[u]^{=} &  0.
		}	
	$$       
	The composition $\tilde{\eta}:=(E_{\eta}\xr{} V)\circ \psi:E_r\xr{} V$ yields a lifting of $\eta$ and proves the exactness of \eqref{equation-representation-Er-easy-sequence}.
       	
	As a consequence, we conclude that $\Hom(E_r,V)$ is a free and finitely generated $R$-module. Furthermore, the morphism
	$
	\Hom(E_r,V)\otimes E_r \xr{} V 
	$	
	is an admissible epimorphism, where $\Hom(E_r,V)\otimes E_r $ is the obvious direct sum of $E_r$ terms. That is, we have an exact sequence in $\mc{U}$:
      	$$
	0\xr{} K \xr{} \Hom(E_r,V)\otimes E_r \xr{} V \xr{} 0.
	$$      

	Let the following be an exact sequence in  $\mc{U}$:
	$$
	0\xr{} X \xr{} Y \xr{} Z \xr{} 0.
	$$
	We will show by induction on $r$ that, for all $s\geq r$ and $Y\in \mc{U}_s,Z\in \mc{U}_r$, the map 
	$\Hom(E_s,Y)\xr{} \Hom(E_s,Z)$ is surjective. Suppose first that $r=0$. Since 
	$$
	\Hom(E_s,Y)\otimes E_s\xr{} Y \xr{} Z
	$$
	is an admissible epimorphism and factors through $\Hom(E_s,Y)\otimes E_0$, the claim follows. 
	
	For $r\geq 1$, we consider:
       	$$
	\xymatrix
	{
		0 \ar[r] &   X \ar[r]\ar[d]^{=}& Y' \ar[r]\ar[d] &  Z_{r-1} \aril[d]\ar[r] & 0 \\
		0 \ar[r] &   X \ar[r]& Y \ar[r]\ar[dr] &  Z \ar[r]\are[d] & 0 \\
		& & & Z_0.
	}
	$$       
	Since the $r=0$ case is proved, we only need to lift morphisms contained in $\Hom(E_s,Z_{r-1})$ to $Y$. We know $Y'\in \mc{U}_t$ for some $t$ and $\Hom(E_t,Y')\xr{} \Hom(E_t,Z_{r-1})$ is surjective by induction. Since every morphism $E_t\xr{} Y$ factors through $E_r$, we are done with the first part of the proposition. The second part follows easily. 
	\end{proof}
\end{proposition}

\begin{proposition}\label{proposition-modules}
  Assumptions and $E_r$ as in Proposition \ref{proposition-representation}. We set $S_r=\Hom(E_r,E_r)^{{\rm op}}$ and $I:=\ker(S_r\xr{} S_0)$. The functor 
 	$$
	\Hom(E_r,-):\mc{U}_r\xr{} \mc{U}_r(\text{($S_r$-modules)},R)
	$$ 
	is an equivalence of categories. Moreover, $I^n=\ker(S_r\xr{} S_{n-1})$, for all $r\geq n\geq 1$, and $I^n/I^{n+1}\cong \left({\rm Ext}^1(\mathbf{1},\mbf{1})^{\vee}\right)^{\otimes n}$.
	\begin{proof}
		Since 	
	  	$
		\Hom(E_r,V)\otimes E_r \xr{} V 
		$
	is an admissible epimorphism for all $V\in \mc{U}_r$, the functor is faithful. In order to show that the functor is full, let $f\in \Hom_{S_r}(\Hom(E_r,V),\Hom(E_r,W))$. We need to prove that it induces a morphism of exact sequences
	$$
	\xymatrix
	{	
	  0 \ar[r] & K_V \ar[r] & \Hom(E_r,V)\otimes E_r \ar[r]\ar[d]^{f\otimes {\rm id}_{E_r}} & V \ar[r] & 0 \\
	  0 \ar[r] & K_W \ar[r] & \Hom(E_r,W)\otimes E_r \ar[r] & W \ar[r] & 0.
	}
	$$	
	Again, we have an admissible epimorphism $\Hom(E_s,K_V)\otimes E_s\xr{} K_V$, for some $s\geq r$. Therefore, it suffices to show that for every $a\in \Hom(E_r,V)\otimes_R \Hom(E_s,E_r)=\Hom(E_r,V)\otimes_R \Hom(E_r,E_r)$ with trivial image in $\Hom(E_r,V)$ under the composition map, $(f\otimes {\rm id}_{\Hom(E_r,E_r)})(a)$ has also trivial image in $\Hom(E_r,W)$. This follows immediately, because $f$ is a morphism of $S_r$-modules.
	
	Finally, let us prove the essential surjectivity. It is sufficient to show that the map
	$$
	\Ext^1(\mbf{1},\mbf{1})\xr{} \Ext_{S_r}(R,R),
	$$
	is an isomorphism for $r\geq 1$. The claim is easily proved for $r=1$. We would like to prove $\ker(S_r\xr{} S_1)=I^2$ for $r\geq 2$, which implies $ \Ext_{S_r}(R,R)=\Ext_{S_1}(R,R)$ and reduces the assertion to the $r=1$ case.  Let us show $\Hom(E_r,T^{\otimes r})=I^r$. Let $X$ be the following pullback
	$$
	\xymatrix{
	  0 \ar[r] & T^{\otimes r} \ar[r]\ar[d]^{=} & X \ar[d] \ar[r] & T^{\otimes r-1} \ar[r]\ar[d] & 0 \\
	  0 \ar[r] & T^{\otimes r} \ar[r] & E_{r} \ar[r]  & E_{r-1} \ar[r]  &  0, 
	}
	$$
	and set $J:=\Hom(E_r,X)\subset S_r$. Since $X=\Gamma^{\otimes r-1}\otimes E_1$, we conclude $I\cdot J= \Hom(E_r,T^{\otimes r})$. By induction, $J$ maps to $I^{r-1}$ in $S_{r-1}$, 
	hence $I^{r}=\Hom(E_r,T^{\otimes r})$.  
	\end{proof}
\end{proposition}

\subsection{} For a fibre functor $\omega:\mc{U}\xr{} \text{($R$-modules)}$, we have 
$$
\End(\omega)=\varprojlim_{r\geq 0} \Hom(\omega_{\mid \mc{U}_r},\omega_{\mid \mc{U}_r}),
$$
where $\Hom(\omega_{\mid \mc{U}_r},\omega_{\mid \mc{U}_r})$ are the $R$-linear natural transformations; we set 
	$$
	I_{\omega}:=\ker(\End(\omega)  \xr{} \End(\omega_{\mid \mc{U}_0})).
	$$

\begin{corollary}\label{corollary-Tanaka}
	Suppose $\mc{U}$ is free. A fibre functor $\omega:\mc{U}\xr{} \text{($R$-modules)}$ induces an equivalence of categories
       	$$
	\mc{U}\xr{} \mc{U}(\End(\omega),R).
	$$       
	Moreover, $I_{\omega}^n=\ker(\End(\omega)  \xr{} \End(\omega_{\mid \mc{U}_{n-1}}))$, and there is a natural isomorphism 
	$I_{\omega}^n/I_{\omega}^{n+1}\cong \left({\rm Ext}^1(\mathbf{1},\mbf{1})^{\vee}\right)^{\otimes n}$.
	\begin{proof}
		The equivalence of categories is an immediate consequence of Proposition \ref{proposition-modules}. The natural isomorphism $I_{\omega}^n/I_{\omega}^{n+1}\cong \left({\rm Ext}^1(\mathbf{1},\mbf{1})^{\vee}\right)^{\otimes n}$ is induced by
		$$\Hom_R(I_{\omega}/I_{\omega}^2,R)\cong \Ext^1_{\End(\omega)}(R,R)=\Ext^1_{\mc{U}}(\mbf{1},\mbf{1}),$$ and the isomorphism $(I_{\omega}/I_{\omega}^{2})^{\otimes n} \xr{\cong} I_{\omega}^n/I_{\omega}^{n+1}$.	
	\end{proof}
\end{corollary}

\begin{definition}\label{definition-F-inv-fibre-functor}(cf.~Definition \cite{BF})
	Let $\sigma:R\xr{} R$ be an endomorphism. Let $F:\mc{U} \xr{} \mc{U}$ be a $\sigma$-linear exact functor, where $\mc{U}$ is a category as in Definition \ref{definition-Ur}. A pair $(\omega,\eta)$, where $\omega:\mc{U}\xr{} \text{($R$-modules)}$  is a fibre functor, and 
	$$
	\eta: \omega\otimes_{R,\sigma} R \xr{\cong} \omega \circ F
	$$
	is an $R$-linear natural isomorphism, is called an $F$-\emph{fibre functor}.
\end{definition}

There is a unique isomorphism $F(\mbf{1})\cong \mbf{1}$ compatible with $\eta$ and $\omega(\mbf{1})\cong R$.

\subsection{} For two $F$-fibre functors $x=(\omega_1,\eta_1),y=(\omega_2,\eta_2)$, we define 
$$
\rho_{y,x}: \Hom_R(\omega_1,\omega_2) \xr{} \Hom_R(\omega_1\otimes_{R,\sigma} R,\omega_2\otimes_{R,\sigma} R) 
$$
to be the $R$-module morphism rendering commutive the following diagram: 
$$
\xymatrix
{
	\omega_1\otimes_{R,\sigma} R \ar[r]^{\eta_1} \ar[d]_{\rho_{y,x}(a)} & \omega_1\circ F \ar[d]^{a\circ F}\\
	\omega_2\otimes_{R,\sigma} R \ar[r]^{\eta_2}  & \omega_2\circ F. 
}
$$

For all $r\geq 0$, we obtain induced maps
$$
\rho_{y,x}: \Hom_R(\omega_{1 \mid \mc{U}_r},\omega_{2 \mid \mc{U}_r}) \xr{} \Hom_R(\omega_{1 \mid \mc{U}_r}\otimes_{R,\sigma} R,\omega_{2\mid \mc{U}_r}\otimes_{R,\sigma} R) 
$$

\begin{remark}
	If $\mc{U}$ is free then 
	$$
	\varprojlim_{r\geq 0} \Hom_R(\omega_{1 \mid \mc{U}_r},\omega_{2 \mid \mc{U}_r}) \otimes_{R,\sigma} R \xr{} \Hom_R(\omega_1\otimes_{R,\sigma} R,\omega_2\otimes_{R,\sigma} R), 
	$$ 
	is an isomorphism, because both modules can be identified with $\varprojlim_{r\geq 0} \omega_2(E_r)\otimes_{R,\sigma} R$, where $(E_r,e_r)$ is representing $\omega_{1\mid \mc{U}_r}$ as in Proposition \ref{proposition-representation}. 
\end{remark}

\subsubsection{} If $x=y$, then $\rho_x:=\rho_{x,x}$ is a morphism of $R$-algebras. If $\omega_1=\omega_2$, then 
\begin{equation}\label{equation-rho-omega}
\rho_{y,x}= (\eta_2^{-1}\circ \eta_1)\circ \rho_{x}.
\end{equation}

We can consider $\Hom(\omega_1,\omega_2)$ as an $\End(\omega_2)$-module. Since there exists a natural isomorphism $\omega_2\xr{} \omega_1$, it is a free $\End(\omega_2)$-module of rank $1$. Obviously, 
\begin{equation}\label{equation-rho-yx}
\rho_{y,x}(ab)=\rho_y(a)\cdot \rho_{y,x}(b)
\end{equation}
for all $a\in \End(\omega_2)$ and $b\in \Hom(\omega_1,\omega_2)$.

\begin{remark}
	Suppose $\mc{U}$ is free, and let $(\omega,\eta)$ be an $F$-fibre functor. Let $(E_r,e_r)$ be representing $\omega_{\mid \mc{U}_r}$, we obtain a projective system $E=(E_r)_r$ with $e_{r+1}\mapsto e_r$, and an isomorphism 
	$$
	\tau:\End(\omega) \xr{\cong}\varprojlim_{r\geq 0} \omega(E_r), \qquad 1\mapsto (e_r)_r=:e.
	$$
	Let $\phi:E\xr{} F(E)$ be the unique morphism of projective systems in $\mc{U}$ such that $\omega(\phi)(e)=\eta(e\otimes 1)$. Via the isomorphism $\tau$, we have 
	$$
	\rho_{(\omega,\eta)}=\varprojlim_{r\geq 0} \left[ \omega(E)\xr{\omega(\phi)} \omega(F(E)) \xr{\eta^{-1}} \omega(E)\otimes_{R,\sigma} R \right].
	$$
\end{remark}

\begin{definition}\label{definition-Hom-F-fibre-functors}
	For two $F$-fibre functors $x=(\omega_1,\eta_1),y=(\omega_2,\eta_2)$, we set
	\begin{align*}
	\Hom(x,y)&:=\{a\in \Hom_R(\omega_1,\omega_2) \mid \rho_{y,x}(a)=a\otimes {\rm id}_R\},\\
	{\rm Isom}(x,y)&:=\{a\in {\rm Isom}_R(\omega_1,\omega_2) \mid \rho_{y,x}(a)=a\otimes {\rm id}_R\}.
	\end{align*}
\end{definition}

\subsubsection{}
Suppose that $\Ext^1(\mbf{1},\mbf{1})$ is free and finitely generated; we set $\Gamma=\Ext^1(\mbf{1},\mbf{1})^{\vee}$. We have a $\sigma$-linear map
	$$
	\Ext^1(F,F):\Ext^1(\mbf{1},\mbf{1}) \xr{} \Ext^1(F(\mbf{1}),F(\mbf{1}))=\Ext^1(\mbf{1},\mbf{1}),  
	$$
	and let $\phi:\Gamma^{\vee}\otimes_{R,\sigma} R \xr{} \Gamma^{\vee}$ be the associated $R$-linear map. Define 
	$$
	\phi^{\vee}:\Gamma \xr{}  (\Gamma^{\vee} \otimes_{R,\sigma} R)^{\vee}=\Gamma\otimes_{R,\sigma} R
	$$
	as dual of $\phi$. 

	\begin{lemma}\label{lemma-rho-formulas}
		Suppose $\mc{U}$ is free. Set $I=\ker(\End(\omega_2)  \xr{} \End(\omega_{2\mid \mc{U}_0}))$ and $H:=\Hom(\omega_1,\omega_2)$.
	\begin{enumerate}
		\item The following map is an isomorphism
			$$
			\Hom_R(\omega_1(\mbf{1}),\omega_2(\mbf{1}))\otimes_R \frac{I^n}{I^{n+1}} \xr{} 	\frac{I^nH}{I^{n+1}H}, \quad b\otimes a\mapsto a\cdot \tilde{b},				
			$$
			where $\tilde{b}\in H$ is a lifting of $b$.
		\item For all $n$, $\rho_{y,x}$ induces a map 
			$$
			\rho_{y,x}^{(n)}:	\frac{I^nH}{I^{n+1}H} \xr{} 	\frac{I^nH}{I^{n+1}H}\otimes_{R,\sigma} R,
			$$
			which equals $\rho_{y,x}^{(0)}\otimes (\phi^{\vee})^{\otimes n}$ via the isomorphism of (1). 	
	\end{enumerate}
	\begin{proof}
		Assertion (1) is obvious. For (2), we may use \eqref{equation-rho-yx} to reduce to the case $x=y=(\omega,\eta)$. Since $\rho_{x}$ is a morphism of $R$-algebras, it suffices to compute for $n=1$. Recall that we have an extension 
			$$
			0\xr{} T \xr{} E_1 \xr{\epsilon} \mbf{1} \xr{} 0.
			$$
			Choose an isomorphism $R\xr{\tau} \omega(\mbf{1})$ and $e\in \omega(E_1)$ such that $\omega(\epsilon)(e)=\tau(1)$. Via $\tau$, we may identify $\omega(T)$, $\Gamma$ and $I/I^2$. For $a\in I/I^2$, $a\cdot e\in \omega(T)$ corresponds to $a$. Moreover, we identify $\mbf{1}$ and $F(\mbf{1})$ via the isomorphism $t:\mbf{1}\xr{\cong} F(\mbf{1})$ satisfying $\omega(t)(\tau(1))=\eta(\tau(1)\otimes 1)$. We have 
			$$ \Hom(\mbf{1},T)\otimes_{R,\sigma} R \xr{} \Hom(F(\mbf{1}),F(T)) \xr{\circ t} \Hom(\mbf{1},F(T)), $$
	 		which yields $F(T)=(\Gamma\otimes_{R,\sigma} R)\otimes \mbf{1}$.		
			
			Let $\psi:E_1\xr{} F(E_1)$ be the unique morphism such that $\omega(\psi)(e)=\eta(e\otimes 1)$. Then 
			$$
			\rho_{x}(a)=\eta^{-1}(\omega(\psi)(a\cdot e))\in \omega(T)\otimes_{R,\sigma} R = I/I^2\otimes_{R,\sigma} R.
			$$
			The restriction of $\psi$ to $T$ is given by the class $f$ of $F(E_1)$ in $\Ext^1(\mbf{1},(\Gamma\otimes_{R,\sigma} R)\otimes \mbf{1})=\Hom_R(\Gamma,\Gamma\otimes_{R,\sigma}R)$. By definition, we have $(\gamma\otimes 1)\circ f=\phi(\gamma\otimes 1)$ for all $\gamma\in \Gamma^{\vee}$, which proves $f=\phi^{\vee}$. 
	\end{proof}
\end{lemma}

\begin{definition}
	Suppose $\mc{U}$ is free. A full subcategory $\mc{U}'$ of $\mc{U}$ is called \emph{admissible} if the following conditions are satisfied:
	\begin{itemize}
		\item $\mbf{1}\in \mc{U}'$, 
		\item there exists a fibre functor $\omega:\mc{U}\xr{} \text{($R$-modules)}$, and a two-sided ideal $K_{\omega}\subset \End(\omega)$, such that
			$$
			E\in \mc{U}' \Leftrightarrow K_{\omega}\cdot \omega(E)=0,
			$$
			and $\End(\omega_{\mid \mc{U}_r})/K_{\omega}\cdot \End(\omega_{\mid \mc{U}_r})$ is an object in $\mc{U}(\End(\omega),R)$ for every $r\geq 0$. 
		\end{itemize}
\end{definition}

We could replace the second condition by requiring the existence of $K_{\omega}$ for every $\omega$. An admissible subcategory is automatically an exact category. 
In view of Corollary \ref{corollary-Tanaka}, the ideal $K_{\omega}\cdot \End(\omega_{\mid \mc{U}_r})$ depends only on $\mc{U}'$, and $\omega$ induces an equivalence of categories 
\begin{equation}\label{equation-equivalence-U'}
\mc{U}'\xr{} \mc{}\mc{U}(\End(\omega)/K_{\omega},R).
\end{equation}
\begin{remark}\label{remark-independence-quotient}
	The $R$-module defined by
$$
\frac{K_{\omega}\cap I^n_{\omega}}{K_{\omega}\cap I^{n+1}_{\omega}} \subset \frac{I^n_{\omega}}{I^{n+1}_{\omega}} = \left(\Ext^1(\mbf{1},\mbf{1})^{\vee}\right)^{\otimes n} 
$$
is independent of the choice of $\omega$. 
\end{remark}

\begin{lemma}\label{lemma-phi-U'}
Suppose $F(\mc{U}')$ is contained in $\mc{U}'$. 
For every $F$-fibre functor $(\omega,\eta)$, and every $n\geq 0$,  
			$$
			\rho_{(\omega,\eta)}(K_{\omega}\cdot \End(\omega_{\mid \mc{U}_n}))\subset K_{\omega}\cdot \End(\omega_{\mid \mc{U}_n})\otimes_{R,\sigma} R,
			$$
			and
			$$
			(\phi^{\vee})^{\otimes n}\left( \frac{K_{\omega}\cap I^n_{\omega}}{K_{\omega}\cap I^{n+1}_{\omega}}  \right) \subset \frac{K_{\omega}\cap I^n_{\omega}}{K_{\omega}\cap I^{n+1}_{\omega}} \otimes_{R,\sigma} R.
			$$
	\begin{proof}
%
		By Corollary \ref{corollary-Tanaka}, there is $E'_n$ in $\mc{U}'\cap \mc{U}_n$ with $\omega(E'_n)\cong \End(\omega_{\mid \mc{U}_n})/K_{\omega}$. For every $a\in K_{\omega}\cdot \End(\omega_{\mid \mc{U}_n})$, we have $\rho_{(\omega,\eta)}(a)\cdot \omega(E'_n)\otimes_{R,\sigma} R=0$, which easily implies $\rho_{(\omega,\eta)}(a)\in K_{\omega}\cdot \End(\omega_{\mid \mc{U}_n}) \otimes_{R,\sigma} R$.

		We know $\rho_{(\omega,\eta)}\left( I^n\cdot  \End(\omega_{\mid \mc{U}_{n}}) \right)\subset I^n\cdot  \End(\omega_{\mid \mc{U}_{n}})\otimes_{R,\sigma} R$. Since $\End(\omega_{\mid \mc{U}_{n}})/K_{\omega}$ is a projective $R$-module for all $n$, we conclude 
	\begin{multline*}
		\left(K_{\omega}\cdot \End(\omega_{\mid \mc{U}_{n}}) \otimes_{R,\sigma} R\right) \cap \left( I^n\cdot  \End(\omega_{\mid \mc{U}_{n}})\otimes_{R,\sigma} R \right) =\\ 
		( K_{\omega}\cap I^n) \cdot \End(\omega_{\mid \mc{U}_{n}})\otimes_{R,\sigma} R,
	\end{multline*}
	which yields the claim by using Lemma \ref{lemma-rho-formulas}.
	\end{proof}
\end{lemma}

\begin{proposition}\label{proposition-criterion-Frob-inv-path-new}
	Suppose $\mc{U}$ is free and $\mc{U}'$ is an admissible $F$-invariant subcategory, that is, $F(\mc{U}')$ is contained in $\mc{U}'$.   
	Let $x=(\omega_1,\eta_1),y=(\omega_2,\eta_2)$ be $F$-fibre functors on $\mc{U}$. 
	If 
	\begin{equation}\label{equation-phi-tensor-power}
				(\phi^{\vee})^{\otimes n}-{\rm id}\otimes 1: \frac{K_{\omega_1}\cap I^n_{\omega_1}}{K_{\omega_1}\cap I^{n+1}_{\omega_1}} \xr{} \frac{K_{\omega_1}\cap I^n_{\omega_1}}{K_{\omega_1}\cap I^{n+1}_{\omega_1}}  \otimes_{R,\sigma} R
	\end{equation}
	is surjective (resp.~injective) for all $n\geq 1$, then
       	$$
	\Hom(x,y) \xr{} \Hom(x_{\mid \mc{U}'},y_{\mid \mc{U}'}) 
	$$       
	is surjective (resp.~injective).
	\begin{proof}
		In the first step we note, for all $n\geq 0$, the surjectivity of 
		\begin{equation}\label{equation-U-to-U'}
		\Hom(\omega_{1\mid \mc{U}_{n}},\omega_{2\mid \mc{U}_{n}}) \xr{} \Hom(\omega_{1\mid \mc{U}'_{n}},\omega_{2\mid \mc{U}'_{n}}). 
	\end{equation}
	Indeed, in order to prove this, we may suppose $\omega_1=\omega_2=\omega$. Since there is $E'_n\in \mc{U}'_n$ with $\omega(E'_n)\cong \End(\omega_{\mid \mc{U}_n})/K_{\omega}$ (as $\End(\omega_{\mid \mc{U}_n})$-module), we may use the equivalence of categories \eqref{equation-equivalence-U'} to conclude. Moreover, we note that the kernels of \eqref{equation-U-to-U'} form a projective system with surjective transition maps. Therefore 
	$$
	\Hom(\omega_{1},\omega_{2}) \xr{} \Hom(\omega_{1 \mid \mc{U}'},\omega_{2\mid \mc{U}'})
	$$
	is surjective.  After lifting an element in $\Hom(x_{\mid \mc{U}'},y_{\mid \mc{U}'})$ to $\Hom(\omega_{1},\omega_{2})$, one can use the surjectivity of \eqref{equation-phi-tensor-power} to alter a lift and obtain an element in $\Hom(x,y)$ (Lemma \ref{lemma-rho-formulas}, Lemma \ref{lemma-phi-U'}). Clearly the injectivity  of \eqref{equation-phi-tensor-power} implies the uniqueness of a lift in $\Hom(x,y)$.
	\end{proof}
\end{proposition}

\section{Unipotent connections for curves}
\label{section-unipotent-connections}

\subsection{}\label{subsection-unipotent-connections}
Let $k$ be a perfect field of characteristic $p$. We denote by $W(k)$ the ring of Witt vectors and by $\sigma$ the Frobenius endomorphism.  Let $\pi:X\xr{} \Spec(W(k))$ be a smooth projective curve such that $W(k)=H^0(X,\OO_X)$, that is, $X$ is geometrically connected. Let $D\subset X$ be a subscheme such that $\pi_{\mid D}$ is finite and \'etale. We denote by $\hat{X}$ and $\hat{D}$ the completions along the special fibres $X_0$ and $D_0$. 

We denote by $\mc{C}$ (resp.~$\hat{\mc{C}}$) the category of locally free coherent $\OO_{X}$-modules (resp.~$\OO_{\hat{X}}$-modules) with logarithmic connection along $D$ (resp.~$\hat{D}$), that is, objects are of the form $(E,\nabla)$ with $\nabla:E\xr{}  E\otimes_{\OO_{X}} \Omega^1_{X/W(k)}(\log D)$  (resp.~$\nabla:E\xr{}  E\otimes_{\OO_{\hat{X}}} \Omega^1_{\hat{X}/W(k)}(\log \hat{D})$) a connection. The categories $\mc{C}$ and $\hat{\mc{C}}$ are exact categories in the evident way, we set $\mbf{1}=(\OO_{X},d)$ (resp.~$\mbf{1}=(\OO_{\hat{X}},d)$), and define $\mc{U}:=\mc{U}(\mc{C},\bf{1})$ (resp.~$\hat{\mc{U}}:=\mc{U}(\hat{\mc{C}},\bf{1})$) with $R=W(k)$. There is an evident completion functor 
$\mc{U}\xr{}\hat{\mc{U}}, E\mapsto \hat{E}.$ We will only work with unipotent logarithmic connections in the following. 

We have  
\begin{multline*}
  \Ext^1_{\mc{U}}(\mbf{1},(E,\nabla))\xr{\cong}H^1(X, E \xr{\nabla}  E\otimes_{\OO_{X}} \Omega^1_{X/W(k)}(\log D)) \\ = H^1(X_0, \hat{E} \xr{\nabla}  \hat{E}\otimes_{\OO_{\hat{X}}} \Omega^1_{\hat{X}/W(k)}(\log \hat{D}))\xr{\cong}\Ext^1_{\hat{\mc{U}}}(\mbf{1},(\hat{E},\nabla)),
\end{multline*}
so that the completion functor is an isomorphism on $\Ext^1(\mbf{1},\mbf{1})$. Moreover, we have an exact sequence
$$
0\xr{} H^0(X,\Omega^1_{X/W(k)}(D)) \xr{} \Ext^1(\mbf{1},\mbf{1}) \xr{} H^1(X,\OO_X) \xr{} 0,
$$
which implies that $\Ext^1(\mbf{1},\mbf{1})$ is a free and finitely generated $R$-module. It is easy to see that there is a functorial injective map 
$$
\Ext^2_{\hat{\mc{U}}}(\mbf{1},(E,\nabla))\inj H^2(X_0, E \xr{\nabla}  E\otimes_{\OO_{\hat{X}}} \Omega^1_{\hat{X}/W(k)}(\log \hat{D})),
$$
and the analogous statement holds for $\mc{U}$. In particular, $\Ext^2(\mbf{1},\mbf{1})=0$ if $D\not= \emptyset$, and we are in the setup of Proposition \ref{proposition-representation}. 

\begin{proposition} \label{proposition-ff}
	Let $U\subset \hat{X}$ be an open, and $(E_1,\nabla),(E_2,\nabla)$ logarithmic (for $\hat{D}$) unipotent connections on $U$. Then
       	$$
	\Hom((E_{1\mid U},\nabla),(E_{2\mid U},\nabla)) \xr{} \Hom((E_{1\mid (U\backslash \hat{D})},\nabla),(E_{2\mid (U\backslash \hat{D})},\nabla)) 
	$$       
	is bijective.
	\begin{proof} 
		Since $E_1^{\vee}\otimes E_2$, with the induced connection, is unipotent again, we may assume $(E_1,\nabla)=(\OO_{\hat{X}},d)$.

		Injective is evident. Hence the question is local; we may suppose $U={\rm Spf}(A)$, $\hat{D}$ consists only of a single point $x_0$, $\hat{D}=V(t)$ for some $t\in A$, and $\Omega_{A/W(k)}=A dt$. Moreover, we may suppose that $E_2$ has an $A$-basis $e_1,\dots,e_r$ such that the connection matrix is strictly sub-diagonal. We denote by $\widehat{A[t^{-1}]}$ the $p$-adic completion. We note that if $f\in \widehat{A[t^{-1}]}$ satisfies $t\cdot \partial_t(f)\in A$, then $f\in A$. Let $\sum_i f_i e_i$ be a horizontal section with $f_i\in \widehat{A[t^{-1}]}$. By induction, we may assume that $f_1,\dots,f_{r-1}\in A$, thus 
		$
		t\cdot \partial_t(f_r) \in A,
		$
		which proves the claim. 
	\end{proof}
\end{proposition}

\subsection{}\label{subsection-Frobenius-pullback} 
Let us define the Frobenius pullback $F:\hat{\mc{U}}\xr{} \hat{\mc{U}}$. Suppose $U\subset \hat{X}\backslash \hat{D}$ is an open, and $\phi,\phi':U\xr{} U$ two lifts of the absolute Frobenius. 
If $(E,\nabla)$ is a unipotent connection on $U$, then it is automatically nilpotent and we obtain a canonical horizontal morphism
\begin{equation}\label{equation-frobenius-to-frobenius}
\phi^*(E)\xr{\cong} \phi'^{*}(E).
\end{equation}
For $(E,\nabla)=(\OO_U,d)$, it is the identity. 

Let $\{U_i\}_{i}$ be an open covering of $\hat{X}$ together with $\phi_i:U_i\xr{} U_i$, a lifting of the Frobenius for each $i$, such that 
\begin{equation}\label{equation-condition-frobenius}
	\phi_i^*(I_{\hat{D} \mid U_i})=I_{\hat{D}\mid U_i}^p,
\end{equation}	
where $I_{\hat{D}}$ is the ideal for $\hat{D}$. Condition \eqref{equation-condition-frobenius} implies that $\phi_i^*(E_{\mid U_i})$ is a \emph{logarithmic} unipotent connection on $U_i$ for every $E\in \hat{\mc{U}}$. In view of \eqref{equation-frobenius-to-frobenius} and Proposition \ref{proposition-ff}, we can glue to a connection $F(E)\in \hat{\mc{U}}$ via the natural morphisms 
$$
\phi_i^*(E_{\mid U_i})_{\mid U_i\cap U_j} \xr{} \phi_j^*(E_{\mid U_j})_{\mid U_i\cap U_j}. 
$$
This construction does not depend (up to natural isomorphisms) on the choice of the covering and the choices for $\phi_i$. We can always find such a covering, because around every point $x_0\in \hat{X}$, there is an open neighborhood  $U$ and an \'etale morphism $f:U\xr{} \hat{\mathbb{A}}^1_{W(k)}$ with $\{x_0\}=f^{-1}(0)$. Any Frobenius lift on $ \hat{\mathbb{A}}^1_{W(k)}$ can be lifted to $U$. 

In this way, we obtain the Frobenius pullback
\begin{equation}\label{equation-Frob-pullback}
  F:\hat{\mc{U}} \xr{} \hat{\mc{U}}. 
\end{equation}

\subsection{} \label{subsection-explicit-Frobenius-to-Frobenius}
In fact, we do not need Proposition  \ref{proposition-ff} in order to define $F$. For $x_0\in D_0$ let $\phi_1,\phi_2$ be two Frobenius lifts on the local ring $\OO_{\hat{X},x_0}$, such that \eqref{equation-condition-frobenius} is satisfied. In the following we will give an explicit description of the natural map 
\begin{equation}\label{equation-Frobenius-to-Frobenius-singular-point}
f:\phi_1^*E_{x_0} \xr{} \phi_2^*E_{x_0},  
\end{equation}
where $E_{x_0}$ is the stalk of a logarithmic unipotent connection defined in a neighborhood of $x_0$. Let $t\in \OO_{\hat{X},x_0}$ be a generator of $I_{\hat{D}}$, we have $\Omega^1_{\hat{X},x_0}= \OO_{\hat{X},x_0} dt$. Let us write $E_{x_0}[t^{-1}]$ for the induced (regular) connection on $\OO_{\hat{X},x_0}[t^{-1}]$. We denote by $\widehat{\OO_{\hat{X},x_0}[t^{-1}]}$ the $p$-adic completion of $\OO_{\hat{X},x_0}[t^{-1}]$. Over $\widehat{\OO_{\hat{X},x_0}[t^{-1}]}$, we know that $f$ is given by Taylor expansions
\begin{align*}
	\widehat{\OO_{\hat{X},x_0}[t^{-1}]}\otimes_{\phi^*_1,\OO_{\hat{X},x_0}[t^{-1}]}E_{x_0}[t^{-1}] &\xr{} \widehat{\OO_{\hat{X},x_0}[t^{-1}]}\otimes_{\phi^*_2,\OO_{\hat{X},x_0}[t^{-1}]}E_{x_0}[t^{-1}] \\
	1\otimes e &\mapsto \sum_{n=0}^{\infty} \frac{(\phi_1^*(t)-\phi_2^*(t))^n}{n!} \otimes \nabla_{\partial_t}^n(e).
\end{align*}
To define this map over $\OO_{\hat{X},x_0}$, we set 
\begin{align}\label{align-formula-explicit-frobenius-to-frobenius}
	\OO_{\hat{X},x_0}\otimes_{\phi^*_1,\OO_{\hat{X},x_0}}E_{x_0} &\xr{} \OO_{\hat{X},x_0}\otimes_{\phi^*_2,\OO_{\hat{X},x_0}}E_{x_0}  \nonumber \\
	1\otimes e &\mapsto \sum_{n=0}^{\infty} \frac{(\frac{\phi_1^*(t)}{\phi_2^*(t)}-1)^n}{n!} \otimes t^n\nabla_{\partial_t}^n(e).
\end{align}
Since $\frac{\phi_1^*(t)}{\phi_2^*(t)}\in \OO_{\hat{X},x_0}^*$ and $\frac{\phi_1^*(t)}{\phi_2^*(t)}-1 \in p\OO_{\hat{X},x_0}$, we conclude $(\frac{\phi_1^*(t)}{\phi_2^*(t)}-1)^n/n! \in  \OO_{\hat{X},x_0}$. Moreover, $t^n\nabla_{\partial_t}^n=P_n(t\nabla_{\partial_t})$ for a polynomial $P_n$ with integral coefficients, thus $t^n\nabla_{\partial_t}^n(e)\in E_{x_0}$. For the convergence of the sum on the right hand side note that $e\in E(U)$ for some open neighborhood $U$ of $x_0$. Suppose that $t,\phi_1,\phi_2$ are defined on $U$, then we can make sense of the right hand side in $\OO_{\hat{X}}(U)\otimes_{\phi^*_2,\OO_{\hat{X}}(U)} E(U)$, because $t^n\nabla_{\partial_t}^n(e)$ converges $p$-adically to $0$. We used an open neighborhood $U$ instead of working with $\OO_{\hat{X},x_0}$, because the latter is not $p$-adically complete.  

In fact, \eqref{align-formula-explicit-frobenius-to-frobenius} equals 
$$
1\otimes e\mapsto \sum_{n=0}^{\infty} \frac{\log(\frac{\phi_1^*(t)}{\phi_2^*(t)})^n}{n!}\otimes (t\nabla_{\partial_t})^n(e),
$$
where $\log(\frac{\phi_1^*(t)}{\phi_2^*(t)})=-\sum_{n=1}^{\infty} \frac{(1-\frac{\phi_1^*(t)}{\phi_2^*(t)})^n}{n}$, which is well-defined and contained in $p\OO_{\hat{X},z_0}$. We will not use this formula and leave the proof to the reader. 

\section{Frobenius fibre functors} \label{section-Frobenius-fibre-functors}
The goal of this section is to construct $F$-fibre functors (see Definition \ref{definition-F-inv-fibre-functor}) for unipotent logarithmic connections on curves. 

\subsection{Points with good reduction}
Recall the setup in Section \ref{subsection-unipotent-connections}. Let $x:{\rm Spf}(W(k))\xr{} \hat{X}$ be a morphism over $W(k)$ such that the topoligical image $x_0$ is not contained in $D_0$ (we call this a point with good reduction). We set $\omega_x(E):=x^*E=E_{x_0}\otimes_{\OO_{\hat{X},x_0},x^*} W(k)$, for every $E\in \hat{\mc{U}}$.  In order to define 
$$
\eta_x: \omega_x(E)\otimes_{ W(k),\sigma} W(k) \xr{} \omega_x(F(E)),
$$
we choose $\{(U_i,\phi_i)\}$ as in Section \ref{subsection-Frobenius-pullback}. For $i$ with $x_0\in U_i$, we have $\phi_i \circ x= y_i \circ \sigma$ for some $y_i:{\rm Spf}(W(k))\xr{} \hat{X}$ with topological image $x_0$, hence an isomorphism
\begin{equation}\label{equation-stupid-eta}
	\omega_{y_i}(E)\otimes_{W(k),\sigma} W(k) \xr{}  \omega_x(\phi^*_i(E)).
\end{equation}
We will simply write $\phi^*_i(E)$ for $\phi^*_i(E_{\mid U_i})$.   
We define $\eta_x$ to be the composition of \eqref{equation-stupid-eta} with the base change (along $\sigma$) of the natural map 
\begin{equation}\label{equation-fibre-to-fibre-same-reduction}
	\omega_x(E) \xr{} \omega_{y_i}(E), 
\end{equation}	
provided by the connection. Recall that we have a functorial map  \eqref{equation-fibre-to-fibre-same-reduction}, because $x$ and $y_i$ are equal modulo $p$, and the connection  is nilpotent. For $j$ with $x_0\in U_j$, there is a commutative diagram
\begin{equation}\label{diagram-def-F-inv-fibre-functor}
\xymatrix
{
	\omega_x(\phi^*_i(E)) \ar[rr]^{\text{\eqref{equation-frobenius-to-frobenius}}} & &
	\omega_x(\phi^*_j(E)) 
	\\
	\omega_{y_i}(E)\otimes_{W(k),\sigma} W(k) \ar[u]
	& &
	\omega_{y_j}(E)\otimes_{W(k),\sigma} W(k) \ar[u]
	\\
	&
	\omega_{x}(E)\otimes_{W(k),\sigma} W(k),\ar[ul] \ar[ur] 
	&	
}
\end{equation}
which implies that $\eta_x$ is well-defined.

\subsection{Points with bad reduction}\label{Section-bad-reduction}
Let $x:{\rm Spf}(W(k))\xr{} \hat{X}$ be a morphism over $W(k)$ such that the topoligical image $x_0$ is contained in $D_0$, but $x$ does not factor through $\hat{D}$. We will essentially use the same method as in the good reduction case, but the maps \eqref{equation-stupid-eta} will depend on the choice of a logarithm for $p$, so that the Frobenius fibre functors obtained in this way are indexed by $(x,a)$ with $a\in p\Z_p$.  

\subsubsection{} We denote by $z:{\rm Spf}(W(k))\xr{} \hat{D}$ the point with topological image $z_0=x_0$, and let $I\subset \OO_{\hat{X},x_0}$ be the associated ideal; set $W:=W(k)$ and $K:={\rm Frac}(W)$. For every generator $t$ of $I$, we have an evident commutative diagram 
$$
\xymatrix{
\varprojlim_n \OO_{\hat{X},x_0}/I^n \ar[r] 
&
\varprojlim_n (\OO_{\hat{X},x_0}/I^n \otimes_W K) 
\\
W[|t|] \ar[u]^{\cong}\ar[r]
&
K[|t|]. \ar[u]^{\rho_t}_{\cong}
}
$$
We define the subring 
$$
A_t:=\{ \sum_{n=0}^{\infty} a_n\frac{t^n}{n!} \mid a_n\in W, \lim_{n\rightarrow \infty} |\frac{a_n}{n!}|r^n=0 \;\text{for all $0<r<1$}\} \subset K[|t|],
$$
and note $W[|t|]\subset A_t$. It is easy to see that $\rho_t(A_t)$ is independent of the choice of $t$, and we denote this ring by $A_z$. The point $x:{\rm Spf}(W(k))\xr{} \hat{X}$ induces an evident ring homomorphism
$$
x^*:A_z\xr{} W.
$$
This also holds for the point $z$. We will need that the kernel of $z^*$ is a divided power ideal.
\begin{lemma}\label{lemma-Az-pd}
	For $f\in A_z$ with $z^*(f)=0$, we have $f^n\in n!\cdot A_z$ for every $n\geq 1$. 
	\begin{proof}
		By induction we may suppose $f^{n-1}/(n-1)!\in A_z$. There is a unique $g\in A_z$ satisfying 
		$$
		dg=\frac{f^{n-1}}{(n-1)!}df, \qquad \qquad z^*(g)=0,
		$$
		hence $f^n=n!\cdot g$.
	\end{proof}
\end{lemma}

Every lift of the Frobenius $\phi^*:\OO_{\hat{X},z_0}\xr{} \OO_{\hat{X},z_0}$ with $\phi^*(I)\subset I$, extends to $A_z$. We have an obvious connection 
\begin{equation}\label{equation-connection-Az}
d:A_z\xr{} A_z \otimes_{\OO_{\hat{X},z_0}} \Omega^1_{\hat{X},z_0}.
\end{equation}
The reason for introducing $A_z$ is that $d$ is surjective and has kernel $W$. In particular, every unipotent connection can be trivialized over $A_z$. Since we work with logarithmic connections we have to go one step further.

\subsubsection{}
We set 
$$
A_{t,\log(t)}:= A_t[\frac{l^n}{n!}\mid n\geq 1]\subset (K\otimes_W A_t)[l],
$$
and extend \eqref{equation-connection-Az} to a logarithmic connection 
\begin{align}\label{align-A-log-d}
	d:A_{t,\log(t)} &\xr{} A_{t,\log(t)}\otimes_{\OO_{\hat{X},z_0}} \Omega^1_{\hat{X},z_0}(\log z) \\
	  d(\frac{l^n}{n!})&=\frac{l^{n-1}}{(n-1)!} \otimes \frac{dt}{t} \qquad \qquad (n\geq 1), \nonumber
\end{align}
thus $l$ plays the role of $\log(t)$.

In order to relate $A_{t,\log(t)}$ and $A_{t',\log(t')}$, we will need the logarithm 
$$
\log:\OO_{\hat{X},z_0}^*\xr{} A_z.
$$ 
For $f\in \OO_{\hat{X},z_0}^*$, we can uniquely write $f=c\cdot g$, with $c\in W^*$ and $z^*(g)=1$. We set 
$$
\log(f)=\log(c)+\log(g),
$$
where $\log(g)\in A_z$ is determined by the properties $z^*(\log(g))=0$ and $d\log(g)=\frac{dg}{g}$. It is given by the series $\log(g)=-\sum_{n=1}^{\infty}\frac{(1-g)^n}{n}$. Note that $x^*\circ \log=\log\circ x^*$ and $\phi^*\circ \log=\log \circ \phi^*$.

Let $t'$ be another generator for $I$, we would like to define
\begin{align*}
	\alpha_{t,t'}: A_{t,\log(t)} &\xr{} A_{t',\log(t')} \\
	\frac{l^n}{n!} &\mapsto \frac{(l+\log(t/t'))^n}{n!}.
\end{align*}
For this, we need to show $\log(t/t')^n \in n!\cdot A_z$ for every $n$. Indeed, since $\log(W^*)\subset pW$, we may suppose $z^*(t/t')=1$ and use Lemma \ref{lemma-Az-pd}. By using the maps $\alpha_{t,t'}$, we glue a ring $A_{z,\log}$ together with the logarithmic connection \eqref{align-A-log-d}.

In order to extend points to $A_{z,\log}$, we need an extension of $\log:W^*\xr{} pW$ to $\log:W\backslash \{0\}\xr{} pW$ such that $\log(xy)=\log(x)+\log(y)$. We can do this by picking a choice $a\in p\Z_p$ for $\log(p)$. 

Our point $x:{\rm Spf}(W(k))\xr{} \hat{X}$ induces 
$$
x^*:A_{t,\log(t)} \xr{} W, \qquad \frac{l^n}{n!}\mapsto \frac{\log(x^*(t))^n}{n!},
$$
for every $t$, and compatibility with $\alpha_{t,t'}$ yields $x^*:A_{z,\log}\xr{} W$.
Here we have used that $x\neq z$. In order to extend $z$, we will need to use tangent vectors, which is done in the next section. 

Every lift of the Frobenius $\phi^*:\OO_{\hat{X},z_0}\xr{} \OO_{\hat{X},z_0}$ with $\phi^*(I)=I^p$, extends naturally to $A_{z,\log}$. Indeed, choose a generator $t$ for $I$, then $\phi(t)=t^pf$ where $f$ is a unit.  In $A_{t,\log(t)}$ we set 
$$
\phi(\frac{l^n}{n!})= \frac{(pl + \log(f))^n}{n!}.
$$
Again, this is well-defined and compatible with $\alpha_{t,t'}$, so that we obtain an induced map on $A_{z,\log}$.

\begin{proposition}
	Let $E$ be a unipotent logarithmic (with respect to $z$) connection over  $\OO_{\hat{X},z_0}$. For the induced (logarithmic) connection $ A_{z,\log}\otimes_{\OO_{\hat{X}},z_0}E$, we denote by       	
	$$
	( A_{z,\log} \otimes_{\OO_{\hat{X},z_0}} E_{z_0})^{\nabla}=\{ \sum_i a_i \otimes e_i \mid \sum_i a_i\nabla(e_i)+ d(a_i)\otimes e_i=0\}
	$$      
	the horizontal sections. The natural map
	$$
	A_{z,\log} \otimes_W (A_{z,\log} \otimes_{\OO_{\hat{X},z_0}} E_{z_0})^{\nabla}	\xr{} A_{z,\log} \otimes_{\OO_{\hat{X},z_0}} E_{z_0} 
	$$
	is an isomorphism.
	\begin{proof}
		First, note that 
		$$
		A_{z,\log}^{\nabla}=A_{z}^{\nabla}=W.
		$$
		Fix a generator $t$ for $I$. Since every (regular) connection has enough horizontal sections over $A_z$, we can find an $A_z$-basis of $A_{z} \otimes_{\OO_{\hat{X},z_0}} E_{z_0}$  such that the connection matrix is of the form $B\frac{dt}{t}$ where $B$ has only entries in $W$ and is strictly sub-diagonal. Working in $A_{t,\log(t)}$, the matrix $\exp(-Bl)$ provides a basis of horizontal sections.   
	\end{proof}
\end{proposition}

This proposition enables us to identify $\omega_x(E)\cong \omega_{x'}(E)$ via the horizontal sections over $A_{z,\log}$, provided that the topological image of $x$ and $x'$ is $z_0$, and $x\neq z,x'\neq z$. Note that this identification depends on the choice of $\log(p)=a$ if $|x^*(t)|\neq |x'^*(t)|$. In case $|x^*(t)| = |x'^*(t)|$, the fibres are canonically identified.

All we need in order to define 
$$
\eta_{x,a}: \omega_x(E)\otimes_{ W(k),\sigma} W(k) \xr{} \omega_x(F(E))
$$
as in the good reduction case is the commutativity of the following diagram, where we have set $E_{z,\log}:=A_{z,\log} \otimes_{\OO_{\hat{X},z_0}} E_{z_0}$,
$$
\xymatrix
{
	A_{z,\log}\otimes_{\OO_{\hat{X},z_0}} \phi_i^*E_{z_0} 
	\ar[rr]^{{\rm id}_{A_{z,\log}}\otimes \text{\eqref{equation-Frobenius-to-Frobenius-singular-point}}} 
	& &
	A_{z,\log}\otimes_{\OO_{\hat{X},z_0}} \phi_j^*E_{z_0} 
	\\
	A_{z,\log} \otimes_{\phi_i^*,A_{z,\log}} E_{z,\log}  \ar[u]^{=} 
      	& &
	A_{z,\log} \otimes_{\phi_j^*,A_{z,\log}} E_{z,\log}  \ar[u]^{=}
	\\
	&
	W\otimes_{\sigma,W} E_{z,\log}^{\nabla} \ar[ru] \ar[lu]
	&
}
$$
The commutativity is not difficult to prove with the help of the explicit description in Section \ref{subsection-explicit-Frobenius-to-Frobenius}.

\subsection{Tangential basepoints} \label{section-tangential-base-points}
Recall the notations from Section \ref{Section-bad-reduction}. Let $z:{\rm Spf}(W(k))\xr{} \hat{D}$ and $\xi\in \Hom_W(I/I^2,W)\backslash \{0\}$. We are going to define an $F$-fibre functor depending on $(z,\xi,a)$, where $\log(p)=a$. If $\xi$ is an isomorphism, then the construction will not depend on $a$. We will give two construction. In the first one, we will use the methods from Section \ref{Section-bad-reduction}. The second construction will follow the general formalism due to Deligne \cite[p.~242]{Deligne} (cf.~\cite[\textsection3]{BF}).  

\subsubsection{} \label{subsection-tangential-base-points-1} We can extend $z^*:A_z\xr{} W$ to $A_{z,\log}$ by defining
$$
(z,\xi)^*: A_{t,\log(t)}\xr{} W, \qquad \frac{l^n}{n!}\mapsto \frac{\log(\xi(t))^n}{n!},
$$
which is compatible with $\alpha_{t,t'}$. This extension depends on $a=\log(p)$ if $\xi$ is not an isomorphism. 

We define $\phi_{T_z}(\xi)\in \Hom_W(I/I^2,W)$ by the property
$$
\phi_{T_z}(\xi)(t)=\sigma^{-1}(\xi(t)^p\cdot z^*(f)), 
$$
where, as usual, $\phi^*$ is a lifting of the absolute Frobenius with $\phi^*(t)=f\cdot t^p$ and $f\in \OO_{\hat{X},z_0}^*$. It is easy to see that $\phi_{T_z}$ does not depend on the choice of $t$. Moreover,
\begin{equation}\label{equation-phi-xi}
\sigma\circ (z,\phi_{T_z}(\xi))^* = (z,\xi)^* \circ \phi^*.
\end{equation}
We set 
$$
\omega_{(z,\xi,a)}(E)=(z,\xi)^*(A_{z,\log}\otimes_{\OO_{\hat{X},z_0}} E_{z_0})=z^*E,
$$
and 
\begin{equation}\label{equation-eta-tangential-base-points}
\eta_{(z,\xi,a)}: \omega_{(z,\xi,a)}(E)\otimes_{W,\sigma} W \xr{} \omega_{(z,\phi_{T_z}(\xi),a)}(E)\otimes_{W,\sigma} W \xr{} \omega_{(z,\xi,a)}(\phi^*(E)),
\end{equation}
where the first map is induced by the trivialization over $A_{z,\log}$, and the second map comes from \eqref{equation-phi-xi}. 
\begin{remark}
	The functor $\omega_{(z,\xi,a)}$ is simply $\omega_z$. However, $\eta_{(z,\xi,a)}$ depends on $(\xi,a)$.  If $\xi$ is an isomorphism then so is $\phi_{T_z}(\xi)$, and $\eta_{(z,\xi,a)}$ does not depend on $a$. We will reprove this fact via the second construction of tangential base points below. 

	The second arrow in \eqref{equation-eta-tangential-base-points} can be identified with the isomorphism $z^*E\otimes_{W,\sigma} W \xr{} z^*\phi^*E$ coming from $\sigma \circ z^*=z^*\circ \phi^*$. For another lift $\phi'$ the diagram 
       	$$
	\xymatrix
	{
		z^* \phi^*E \ar[rr]^{\eqref{equation-Frobenius-to-Frobenius-singular-point}} 
		&&
	        z^* \phi'^*E 
		\\
		&
		z^*E\otimes_{W,\sigma} W \ar[ul] \ar[ur]
		&
	}
	$$       
	does not commute in general. 
\end{remark}

\subsubsection{}\label{section-tangential-Deligne}
Let us explain the second construction. We set $T_z=\Hom_W(I/I^2,W)$, and $\P_z:=\P(T_z\oplus W)={\rm Proj}({\rm Sym}^*_W(I/I^2\oplus W))$. We let $D_z=\{0\}\cup \{\infty\}\subset \P_z$. For any generator $t$ of $I/I^2$ we get a corresponding open immersion
$$
\Spec \; W[t] \xr{} \P_z 
$$
with $\infty$ as complement.
As usual, we denote by $\hat{\P}_z$ the $p$-adic completion. We denote by $\mc{U}_z$ (resp.~$\hat{\mc{U}}_z$) the category of unipotent connections on $\P_z$ (resp.~$\hat{\P}_z$) with logarithmic singularities at $D_z$.

We have a residue functor
\begin{align*}
{\rm Res}_z:\hat{\mc{U}} &\xr{} \text{(free f.g.~$W$-modules equipped with a nilpotent endomorphism)}, \\
(E,\nabla)&\mapsto (z^*E,{\rm Res}_z(\nabla)).
\end{align*}
For $E\in \hat{\mc{U}}$ we obtain a unipotent connection $(z^*E\otimes_W \OO_{W[t]},d+{\rm Res}_z(\nabla)\frac{dt}{t})$ on $\Spec(W[t])$ with logarithmic singularities at $0$. Any two extensions of this connection to an object in $\mc{U}_z$ are canonically isomorphic, because taking residues at $0$ is an equivalence of categories for $\mc{U}_z$. This construction is independent of the choice of $t$, because $\frac{dt}{t}$ is. Therefore, we obtain a functor
$$
\mathcal{T}_z: \hat{\mc{U}} \xr{} \hat{\mc{U}}_z.
$$
Obviously, $\mathcal{T}_z$ is already defined on unipotent logarithmic connections over $\OO_{\hat{X},z_0}$. 

Let $\phi^*$ be as in Section \ref{subsection-tangential-base-points-1}. We get an associated lifting of the Frobenius 
$$
\phi_{T_z}:\Spec \; W[t] \xr{} \Spec \; W[t], \quad \phi_{T_z}^*(t)=z^*(f)t^p, \quad \phi_{T_z\mid W}^*=\sigma,
$$
and can extend it to $\P_z$. There is a unique isomorphism 
\begin{equation}\label{equation-Frob-cont-tang}
	\mathcal{T}_z\phi^*(E)\xr{} \phi_{T_z}^*\mathcal{T}_z(E)
\end{equation}
inducing the identity on the fibre at $0$, when we identify 
$$0^*\phi_{T_z}^*\mathcal{T}_z(E)=0^*\mathcal{T}_z(E)\otimes_{W,\sigma} W=z^*E\otimes_{W,\sigma} W, \quad 0^*\mathcal{T}_z\phi^*(E)=z^*\phi^*E=z^*E\otimes_{W,\sigma} W.$$
In order to see that \eqref{equation-Frob-cont-tang} gives rise to a natural isomorphism 
\begin{equation}\label{equation-Frob-cont-tang-2}
	\mathcal{T}_z\circ F\xr{} F\circ \mathcal{T}_z
\end{equation}
of functors $\hat{\mc{U}} \xr{} \hat{\mc{U}}_z$, we have to show that 
$$
\xymatrix
{
	z^* \phi^*E \ar[d]_{\eqref{equation-Frobenius-to-Frobenius-singular-point}} \ar[rr]^{\text{\eqref{equation-Frob-cont-tang}}}
	& & 
	0^*\phi_{T_z}^*\mathcal{T}_z(E) \ar[d]^{\eqref{equation-Frobenius-to-Frobenius-singular-point}}
	\\
	z^* \phi'^*E  \ar[rr]^{\text{\eqref{equation-Frob-cont-tang}}} 
	& &
	0^*{\phi'_{T_z}}^*\mathcal{T}_z(E) 
}
$$
is commutative. This follows from the explicit description of \eqref{equation-Frobenius-to-Frobenius-singular-point} in Section \ref{subsection-explicit-Frobenius-to-Frobenius}.

Given $\mathcal{T}_z$ and \eqref{equation-Frob-cont-tang-2}, we can make every $F$-fibre functor on $\hat{\mc{U}}_z$ to an $F$-fibre functor on $\mc{U}$. In particular, every $\xi\in \Hom(I/I^2,W)\backslash \{0\}\subset \P_z(W)$ gives rise to an $F$-fibre functor. If $\xi$ is an isomorphism then $\xi$ has good reduction and there is no dependence on $\log(p)$. 

For the compatibility of the first construction in Section \ref{subsection-tangential-base-points-1} and the second, we note that \eqref{equation-eta-tangential-base-points} corresponds to
$$
\eta_{(z,\xi,a)}: z^*E\otimes_{W,\sigma} W \xr{} z^*E\otimes_{W,\sigma} W \xr{} z^*\phi^*E
$$
by using $\omega_{(z,\xi,a)}=\omega_z$, and it is not difficult to show that the first map is given by $\exp((\log(\xi(t))-\log(\phi_{T_z}(\xi)(t))) \cdot {\rm Res}_z(\nabla))\otimes {\rm id}_W$ for any generator $t$ of $I/I^2$.

\section{Proof of the main theorem}\label{section-main-thms}
Recall the setup of Section \ref{subsection-unipotent-connections}. 

\begin{lemma}\label{lemma-Res}
	Let $a\in \hat{D}(W(k))$. Let $x=(\omega_a,\eta)$ be a tangential $F$-fibre functor at the point $a$. We denote by $\hat{\mc{U}}'$ the unipotent logarithmic connections on $\hat{X}$ with respect to $\hat{D}-a$.
	\begin{enumerate}
		\item ${\rm Res}_a(\nabla)$ defines an element in $\End(\omega_a)$ with $\rho_x({\rm Res}_a(\nabla))=p\cdot {\rm Res}_a(\nabla)\otimes \id_{W(k)}$.
		\item If $\hat{D}\neq a$ then ${\rm Res}_a(\nabla)\in I_{\omega_a}\backslash  I_{\omega_a}^2$ and $I_{\omega_a}/ (W(k)\cdot {\rm Res}_a(\nabla) +  I_{\omega_a}^2)$ is $p$-torsion free.
		\item If $\hat{D}=a$ then ${\rm Res}_a(\nabla)\in I^2_{\omega_a}\backslash  I_{\omega_a}^3$ and $I^2_{\omega_a}/ (W(k)\cdot {\rm Res}_a(\nabla) +  I_{\omega_a}^3)$ is $p$-torsion free.
	\end{enumerate}
	In any case, $\hat{\mc{U}}'$ is an admissible subcategory of $\hat{\mc{U}}$ and the corresponding ideal is generated by ${\rm Res}_a(\nabla)$.
	\begin{proof}
		(1) is easy to check. For (2). Since $\hat{D}\neq a$, ${\rm Res}_a:H^0(\hat{X},\omega_{\hat{X}}(\hat{D}))\xr{} W(k)$ is surjective. For $b\in W(k)^*$ choose $\omega\in H^0(\hat{X},\omega_{\hat{X}}(\hat{D}))$ with ${\rm Res}(\omega)=b$. Then $E=\OO_{\hat{X}}e_0\oplus \OO_{\hat{X}}e_1$ with $\nabla(e_0)=\omega \otimes e_1, \nabla(e_1)=0,$ shows the claim. 
		
		For (3). In this case, we have $H^1_{{\rm dR}}(\hat{X}/W)=\Ext^1_{\hat{\mc{U}}'}(\mbf{1},\mbf{1})=\Ext^1_{\hat{\mc{U}}}(\mbf{1},\mbf{1})$. Hence ${\rm Res}_a(\nabla)\in I^2_{\omega_a}$. We have 
		$$
		I^2_{\omega_a}/ I_{\omega_a}^3 =\left(\Ext^1_{\hat{\mc{U}}}(\mbf{1},\mbf{1})^{\otimes 2}\right)^{\vee}=\left(H^1_{{\rm dR}}(\hat{X}/W)^{\otimes 2}\right)^{\vee},
		$$
		and the image of ${\rm Res}_a(\nabla)$ in $I^2_{\omega_a}/ I_{\omega_a}^3$ is given by the Poincar\'e pairing. We are done with (3), because this pairing is perfect for curves, that is, induces an isomorphism
		$$
		H^1_{{\rm dR}}(\hat{X}/W) \cong H^1_{{\rm dR}}(\hat{X}/W)^{\vee}.
		$$

		For the last claim on $\hat{\mc{U}}'$, we can take $K_{\omega_a}$ to be the two sided ideal generated by ${\rm Res}_a(\nabla)$. Since the leading term of ${\rm Res}_a(\nabla)$ is not divisible by $p$, we obtain that $I^n/(K_{\omega}\cap I^n + I^{n+1})$ is a free $W(k)$-module for all $n$. By induction on $n$, we see that $\End(\omega_{a\mid \mc{U}_n})/K_{\omega}$ is contained in $\mc{U}(\End(\omega_a),W(k))$.
	\end{proof}
\end{lemma}

\begin{thm}\label{thm-reduction-to-projective-curve}
	Let $k$ be a perfect field of characteristic $p$. Suppose $X/W(k)$ is a smooth projective geometrically connected curve. Let $D=\sum_i a_i$, with $a_i\in X(W(k))$, be a divisor that is \'etale over $W(k)$. Let $x,y$ be two $F$-fibre functors (Definition \ref{definition-F-inv-fibre-functor}) on the category $\hat{\mc{U}}$ of unipotent connections on $\hat{X}$ with logarithmic singularities at $\hat{D}$. We denote by $\hat{\mc{U}}'$ the category of  unipotent connections on $\hat{X}$. Then
$$
\Hom(x,y)\xr{} \Hom(x_{\mid \hat{\mc{U}}'},y_{\mid \hat{\mc{U}}'})
$$
is bijective (see Definition \ref{definition-Hom-F-fibre-functors}).
\begin{proof}
	We denote by $\hat{\mc{U}}_{a_1}$ the category of unipotent connections on $\hat{X}$ with logarithmic singularities at $\hat{D}-a_1$. It suffices to prove that 
	$$
	\Hom(x,y)\xr{} \Hom(x_{\mid \hat{\mc{U}}_{a_1}},y_{\mid \hat{\mc{U}}_{a_1}})
	$$
	is bijective. Since $\hat{\mc{U}}_{a_1}$ is an admissible subcategory of $\hat{\mc{U}}$ by Lemma \ref{lemma-Res}, we may use Proposition \ref{proposition-criterion-Frob-inv-path-new} and show that \eqref{equation-phi-tensor-power} is bijective. In view of Remark \ref{remark-independence-quotient} and Lemma \ref{lemma-rho-formulas}(2) applied to a tangential fibre functor at $a_1$, we conclude from Lemma \ref{lemma-Res}(1) that 
$$
(\phi^{\vee})^{\otimes n}\left( \frac{K_{\omega_x}\cap I^n_{\omega_x}}{K_{\omega_x}\cap I^{n+1}_{\omega_x}} \right) \subset p\cdot \left(\frac{K_{\omega_x}\cap I^n_{\omega_x}}{K_{\omega_x}\cap I^{n+1}_{\omega_x}}  \otimes_{W(k),\sigma} W(k)\right).
$$
Hence \eqref{equation-phi-tensor-power} is bijective. 
\end{proof}
\end{thm}

\begin{thm}\label{thm-construction-path}
	Suppose $X/W(k)$ is either $\P^1_{W(k)}$ or an elliptic curve.   Let $D=\sum_i a_i$, with $a_i\in X(W(k))$, be a divisor that is \'etale over $W(k)$. Let $x,y$ be two $F$-fibre functors  on the category $\hat{\mc{U}}$ of unipotent connections on $\hat{X}$ with logarithmic singularities at $\hat{D}$ as constructed in Section \ref{section-Frobenius-fibre-functors}, that is, $x$ (resp.~$y$) is attached to a point with good or bad reduction or is tangential. Let $x_0\in X(k)$ (resp.~$y_0\in X(k)$) be the reduction of the underlying point of $x$ (resp.~$y$). If $X$ is an elliptic curve then we assume that $y_0-x_0$ has order prime to $p$. 
	
	Then there exists $\gamma_{y,x}\in {\rm Isom}(x,y)$ such that $\gamma_{y,x}$ induces the identity on $(\OO_{\hat{X}},d)$.  
	\begin{proof}
		By using Theorem \ref{thm-reduction-to-projective-curve}, we may suppose $D=\emptyset$. 

		Let $x=(\omega_a,\eta_a), y=(\omega_b,\eta_b)$ with $a,b\in X(W(k))$ the underlying points. We have identified $\omega_a(\OO_{\hat{X}})$ with $\omega_b(\OO_{\hat{X}})$ by $1\mapsto 1$. If $X=\P^1_{W(k)}$ then there is nothing to prove, because every unipotent connection is a direct sum of trivial ones. 
		
		Now assume that $X$ is an elliptic curve. Certainly we may suppose that $a$ is the zero. For any $b'\in X(W(k))$ in the same residue disc with $b$, there is a canonical isomorphism in ${\rm Isom}((\omega_{b'},\eta_{b'}),(\omega_{b},\eta_{b}))$ inducing the identity on $(\OO_{\hat{X}},d)$, because the fibres are identified via the connection. By Hensel's lemma, there is a prime-to-$p$ torsion point $b'\in X(W(k))$ in the residue disc of $b$. We may assume $b'=b$. 

		Let $N\geq 1$ be such that $(N,p)=1$ and $N\cdot b=0$. Multiplication by $N$, $m_N:\hat{X}\xr{} \hat{X},$ induces a functor 
	$$
	m_N^*: \text{(unipotent connections on $\hat{X}$)} \xr{} \text{(unipotent connections on $\hat{X}$)}. 
	$$	
	This functor is an equivalence of categories, because on $\Ext^1(\mbf{1},\mbf{1})$ and $\Ext^2(\mbf{1},\mbf{1})$ it induces multiplication by $N$ and $N^2$, respectively. Moreover there is a natural isomorphism 
	$$
	F\circ m_N^* \xr{} m_N^* \circ F.
	$$
	Indeed, if $\phi$ is a lifting of the absolute Frobenius on an open $U$ of $\hat{X}$ and $\phi'$ is a lifting of the Frobenius on $m_{N}^{-1}(U)$, then there is a natural isomorphism 
      	$$
	\phi'^*\circ m_N^* \xr{\cong} m_N^* \circ \phi^* 
	$$      
	provided by the connection, because $m_N\circ \phi'$ and $\phi \circ m_N$ agree on the special fibre. 
	
	For a unipotent connection $E$, find $E'$ and $m_N^*E'\cong E$. We obtain a map
      	$$
		\omega_0(E) \xr{\cong} \omega_0(m_N^*E')= \omega_0(E')= \omega_b(m_N^* E')\xr{\cong} \omega_b(E),
	$$      
	which is independent of the choices made. This construction yields the desired element in ${\rm Isom}((\omega_{b'},\eta_{b'}),(\omega_{b},\eta_{b}))$. 
	\end{proof}	
\end{thm}

\begin{remark}\label{remark-same-residue-disc}
  \begin{enumerate}
\item If $X=\P^1_{W(k)}$ then $\gamma_{y,x}$ is unique. In particular, this yields compatibility with composition:
       	$$
	\gamma_{z,x}=\gamma_{z,y}\circ \gamma_{y,x}.
	$$       
\item If $x_0=y_0\not\in D(k)$ then $\gamma_{y,x}$ is simply the usual isomorphism of fibres induced by the connection.
\item Suppose $x_0=y_0\in D(k)$, hence $x$ (resp.~$y$) is attached to a point with bad reduction or is tangential. Let $z\in D(W(k))$ be such that $z_0=x_0=y_0$. Then $\gamma_{y,x}$ fits into a commutative diagram
 $$
 \xymatrix{
 \omega_{x}(E) \ar[rr]^{\gamma_{y,x}}
 &&
 \omega_{y}(E)\\
 &
 E_{z,\log}^{\nabla},\ar[ul] \ar[ur] 
& 
}
$$ 
where $E_{z,\log}^{\nabla}$ are the horizontal sections of $E_{z_0}$ over the ring $A_{z,\log}$ (see Section \ref{Section-bad-reduction}). 
\end{enumerate}
\end{remark}

For the proof of this remark we can reduce to $D=\emptyset$ by Theorem \ref{thm-reduction-to-projective-curve} again. Then the claims follow immediately from the construction of $\gamma_{y,x}$. 

\subsection{Compatibilities}
In this section we sketch the compatibility of the path constructed in Theorem \ref{thm-construction-path} with Besser's Frobenius invariant path.  Again, we assume the setup of Section \ref{subsection-unipotent-connections}, and suppose that $D=\sum_i a_i\neq 0$, with $a_i\in X(W(k))$, is \'etale over $W(k)$. 

We set $K:={\rm Frac}(W(k))$. For a $W(k)$-linear category $\mc{C}$, we denote by $\mc{C}\otimes \Q$ the $K$-linear category with the same objects as in $\mc{C}$ and 
$$
\Hom_{\mc{C}\otimes \Q}(Z,Y):=\Hom_{\mc{C}}(Z,Y)\otimes_{\Z} \Q. 
$$
\subsubsection{}
We set $A:=\OO_{X}(X\backslash D)$. Let us denote by $\hat{A}$ the $p$-adic completion and by $A^{\dagger}\subset \hat{A}$ the weak completion \cite[\textsection1]{MW}. We denote by $\mc{U}(A^{\dagger}\otimes K)$ (resp.~$\mc{U}(\hat{A}\otimes K)$) the unipotent $A^{\dagger}\otimes K$-connections (resp.~$\hat{A}\otimes K$-connections). We have a commutative diagram of functors
$$
\xymatrix
{
	&
	\mc{U}\otimes \Q \ar[rd]^{E\mapsto \hat{E}} \ar[ld]_{E\mapsto E^{\dagger}} 
	&
	\\
	\mc{U}(A^{\dagger}\otimes K)\ar[dr]_{J^{\dagger}}
	&
	&
	\hat{\mc{U}}\otimes \Q \ar[dl]^{\hat{J}}
	\\
	&
	\mc{U}(\hat{A}\otimes K),
	&
}
$$
where $J^{\dagger}$ and $\hat{J}$ are the evident restriction functors. The functors $E\mapsto \hat{E}$ and $E\mapsto E^{\dagger}$ are equivalences of categories. The functors $\hat{J}$ and $J^{\dagger}$ are fully faithful by Proposition \ref{proposition-ff}.

We have a functor $F: \mc{U}(A^{\dagger}\otimes K) \xr{} \mc{U}(A^{\dagger}\otimes K)$, which is defined by $E\mapsto \phi^*E$ if $\phi:A^{\dagger}\xr{} A^{\dagger}$ is a lifting of the absolute Frobenius. For all $E\in \mc{U}$, we may identify $J^{\dagger}(F(E^{\dagger}))$ and $\hat{J}(F(\hat{E}))$. 

\subsubsection{}
Suppose $k$ is a field with $p^n$ elements. Besser \cite[p.~26]{B} constructs for all $x_0\in X(k)$ a fibre functor (in the Tannakian sense)
$$
\omega^B_{x_0}: \mc{U}(A^{\dagger}\otimes K) \xr{} \text{($K$-vector spaces)}
$$
equipped with an evident isomorphism $\omega^B_{x_0} \circ F^{n}\cong \mbf{\omega}^B_{x_0}$. For every pair $x_0,y_0\in X(k)$, there is a unique natural isomorphism $\gamma_{y_0,x_0}^B\in {\rm Isom}^{\otimes}(\omega^B_{x_0},\omega^B_{y_0})$  such that 
$$
\xymatrix
{
	\omega^B_{x_0} \ar[r]^-{\cong} \ar[d]_{\gamma^B_{y_0,x_0}}
	&
	\omega^B_{x_0}\circ F^{n} \ar[d]^{\gamma^B_{y_0,x_0}\circ F^n}
	\\
	\omega^B_{y_0} \ar[r]^-{\cong} 
	&
	\omega^B_{y_0} \circ F^{n} 	
}
$$
is commutative; it is called the \emph{Frobenius invariant path}. For $E\in \mc{U}$ and $x_0\in X(k)$, the following holds.
\begin{description}
	\item[(Good reduction)] If $x_0\in X(k)\backslash D(k)$ then 
		$$\omega^B_{x_0}(E)=\{s\in \hat{E}_{x_0}(]x_0[) \mid \nabla(s)=0\}.$$ 
	Therefore  Remark \ref{remark-same-residue-disc}(2) implies the commutivity of the diagram
	\begin{equation}\label{diagram-Besser-comparison}
	\xymatrix{
		& 
		\omega^B_{x_0}(E) \ar[dl] \ar[dr]
		&
		\\
		\omega_{x}(E)\otimes K \ar[rr]^{\gamma_{x',x}}
		&
		&
		\omega_{x'}(E)\otimes K 
	}
	\end{equation}	
	for all $x,x'\in X(W(k))$ reducing to $x_0$. Here we have used the notation $\hat{E}_{x_0}(]x_0[)=A(]x_0[)\otimes_{\OO_{\hat{X},x_0}} \hat{E}_{x_0}$ with 
	$$
	A(]x_0[)=\{ \sum_{i=0}^{\infty}  a_it^i \mid \text{$a_i\in K$ and $\lim_{i\rightarrow \infty} |a_ir^i| =0$ for all $0<r<1$}\},
	$$
	where $t\in \OO_{\hat{X},x_0}$ reduces to a uniformizing element modulo $p$.
\item[(Bad reduction/Tangential)] If $x_0\in D(k)$ then  $\omega^B_{x_0}(E)=\hat{E}_{z,\log}^{\nabla}\otimes K$, where $z\in D(W(k))$ is the lifting of $x_0$ \cite[p.~41]{B}. Remark \ref{remark-same-residue-disc}(3) implies the commutativity of \eqref{diagram-Besser-comparison} if $x$ (resp.~$x'$) is an $F$-fibre functor attached to a point with reduction $x_0$ or is tangential.
\end{description}
Therefore we obtain an isomorphism
\begin{equation}\label{equation-comparison-with-Besser}
{\rm Isom}(\omega_x\otimes K,\omega_y\otimes K) \xr{\cong} {\rm Isom}(\omega_{x_0}^B,\omega_{y_0}^B),
\end{equation}
for all $F$-fibre functors $x,y$ constructed in Section \ref{section-Frobenius-fibre-functors}, and where $x_0\in X(k)$ denotes the reduction of the point underlying $x$.  

\subsubsection{} For every $F$-fibre functor $x=(\omega,\eta)$ (on $\hat{\mc{U}}$) we obtain an $F^n$-fibre functor $x^{[n]}=(\omega,\eta^{[n]})$, where
$$
\eta^{[n]}=(\omega \otimes_{W(k),\sigma^n} W(k) \xr{\eta} (\omega\circ F) \otimes_{W(k),\sigma^{n-1}} W(k) \xr{\eta^{[n-1]}\circ F} \omega\circ F^n).
$$
We denote by $x^{[n]}\otimes K=(\omega\otimes K, \eta^{[n]}\otimes K)$ the induced $F^n$-fibre functor on $\hat{\mc{U}}\otimes \Q$. 

In fact, \eqref{equation-comparison-with-Besser} induces 
\begin{equation} \label{equation-comparison-with-Besser-2}
	{\rm Isom}(x^{[n]}\otimes K, y^{[n]}\otimes K) \xr{\cong} {\rm Isom}((\omega_{x_0}^B,\omega_{x_0}^B\xr{\cong} \omega_{x_0}^B\circ F^n),(\omega_{y_0}^B,\omega_{y_0}^B\xr{\cong} \omega_{y_0}^B\circ F^n)),
\end{equation}
where ${\rm Isom}$ is as in Definition \ref{definition-Hom-F-fibre-functors}, on the left hand side for $\hat{\mc{U}}\otimes \Q$, and on the right hand side for $\mc{U}(A^{\dagger}\otimes K)$. Although we have explained the construction for $D\neq \emptyset$ only, we also have \eqref{equation-comparison-with-Besser-2} in the case $D=\emptyset$. 

\begin{proposition}\label{proposition-comparison}
	Via \eqref{equation-comparison-with-Besser-2}, the element $\gamma_{y,x}$ constructed in Theorem \ref{thm-construction-path} maps to Besser's Frobenius invariant path.
	\begin{proof}
			Indeed, for $D=\emptyset$, we can see from the explicit description of $\gamma_{y,x}$ that it is a $\otimes$-isomorphism, hence maps to the Frobenius invariant path by uniqueness of the latter.

			For $D\neq \emptyset$, the Frobenius invariant path provides an extension of the restriction of $\gamma_{y,x}$ to $\hat{\mc{U}}'\otimes \Q$, where $\mc{U}'$ is the category of unipotent connections on $X$. By a variant of Theorem \ref{thm-reduction-to-projective-curve} for $\hat{\mc{U}}\otimes \Q$ and $\hat{\mc{U}}'\otimes \Q$, whose proof is essentially the same, this extension is unique. 	
	\end{proof}
\end{proposition}

\section{Application to $p$-adic multiple zeta values}\label{section-application-mutliple-zeta-values}

\subsection{}
Let $X=\P^1_R$, where $R=W(k)$, and $k$ is a finite field of characteristic $p$. Let $D=\coprod_{i=0}^s a_i$, with $a_i\in \P^1(R)$, be a non-empty subscheme of $\P^1_R$; in particular, we have $\bar{a}_i\neq \bar{a}_j$, for $i\neq j$, and where $\bar{a}_i$ denotes the induced $k$-valued point.  

We define $\omega_i\in H^0(\P^1_{R},\Omega^1(D))$, for $i=1,\dots,s$, to be the form with residue $-1$ at $a_i$, and vanishing residues at all other $a_j$ with $j>0$. If $a_0=\infty$ then $\omega_i=\frac{dx}{a_i-x}$. We set $\Gamma=H^0(\hat{\P}^1_{R},\Omega^1(D))^{\vee}=H^1(\hat{\P}^1_{R},\Omega^{*}(\log D))^{\vee}$, $\Delta_r=\bigoplus_{i=0}^r \Gamma^{\otimes i}$, and we write $\omega^{\vee}_{1},\dots,\omega^{\vee}_{s}$ for the basis dual to $\omega_1,\dots,\omega_s$. We equip
$$
E_r:=\OO_{\hat{\P}^1_R}\otimes_{R} \Delta_r
$$
with the connection induced by
$$
\nabla(1\otimes \alpha):= \sum_{i=1}^s \omega_i\otimes \omega_i^{\vee}\cdot \alpha
$$
for all $\alpha\in \Delta_r$. It is a unipotent logarithmic connection and satisfies the properties of Proposition \ref{proposition-representation}. 

\subsection{} Fix a non-degenerate tangent vector $\xi$ at $a_1$, that is $\xi\in {\rm Isom}_R(I_{a_1}/I_{a_1}^2,R)$, where $I_{a_1}$ is the ideal for $a_1$. Let $y\in \P^1(R)$ be such that the reduction $\bar{y}\in \P^{1}(k)$ is not contained in $D(k)$. We identify $\omega_{a_1}(E_r)$ and $\omega_y(E_r)$ with $\Delta_r$ in the evident way.
Note that 
$$
\gamma_{y,(a_1,\xi),E_r}(1)=1+\sum_{i_1,\dots,i_j} (-1)^j \int_{(a_1,\xi)}^y \omega_{i_1}\circ \dots \circ \omega_{i_j} \cdot \omega_{i_1}^{\vee}\otimes \cdots \otimes \omega_{i_j}^{\vee}, 
$$ 
where  $\int_{(a_1,\xi)}^y \omega_{i_1}\circ \dots \circ \omega_{i_j}$ is the iterated Coleman integral (Proposition \ref{proposition-comparison}).

\newcommand{\fps}[1]{\langle \langle #1 \rangle \rangle} 

We denote the prosystem of connections $(E_r,\nabla)$ by $(E,\nabla)$, so that 
$
\gamma_{y,(a_1,\xi),E}
$ 
is an endomorphism of $R\fps{\Gamma}=\varprojlim_r \Delta_r$.
 
\begin{thm}\label{thm-main}
Let $p^{[r]}\subset W(k)$ be the ideal generated by all $\frac{p^i}{i!}$ with $i\geq r$. We have
\begin{equation}\label{equation-to-prove-in-main-thm}
\gamma_{y,(a_1,\xi),E}(1)\in \prod_{r\geq 0} p^{[r]} \cdot \Gamma^{\otimes r}. 
\end{equation}
\begin{proof}
Let us choose a covering $\hat{\P}^1=\bigcup_{i=0}^s U_i$ with open sets such that $\bar{a}_j\not\in U_i$ for all $i\neq j$. Furthermore, we may assume that $U_1$ is the $p$-adic completion of $\P_{W(k)}^1\backslash \{a_0,a_2,\dots,a_s\}$. We can find on each $U_i$ a lifting of the absolute Frobenius $\phi_i$ with the property $\phi_i^*(I_{a_i})=I_{a_i}^p$. 

The construction of $F^*(E,\nabla)$ works as follows. We glue $\phi_i^*E_{\mid U_i}$ via the canonical morphisms
\begin{align*}
g_{ji}:  \phi_i^*&E_{\mid U_i\cap U_j} \xr{}  \phi_j^*E_{\mid U_j\cap U_i} \\
        \phi_i^*&(e)\mapsto \sum_{m=0}^{\infty} \frac{(\phi^*_i(t)-\phi^*_j(t))^{m}}{m!}\cdot \phi_j^*(\nabla_{\partial_t}^m(e)),
\end{align*}
where $t$ is a coordinate on $U_i\cap U_j$ and $e$ is any section of $E$. Note that $\phi_i(t)-\phi_j(t)\in p\cdot \OO(U_i\cap U_j)$ so that the fraction makes sense. Convergence of the sum is evident if $p\neq 2$, and follows from the quasi-nilpotence of the connection if $p=2$. 

Let us identify 
\begin{equation*}
\epsilon_i:\OO_{U_i}\hat{\otimes}_{\sigma,W(k)} W(k)\fps{\Gamma} \xr{\cong} \phi_i^*E_{\mid U_i}, \qquad
\epsilon_i(1\otimes \alpha) = \phi_i^*(\alpha), 
\end{equation*}
for every $\alpha\in W(k)\fps{\Gamma}$.
We claim that
$$
F_{i}=\prod_{r\geq 0}  \OO_{U_i}\otimes_{\sigma, W(k)} p^{[r]}\cdot \Gamma^{\otimes r} \overset{\epsilon_i}{\hookrightarrow} \phi_i^*E_{\mid U_i} 
$$
inherits the connection from $\phi_i^*E_{\mid U_i}$ and glues via $\{g_{ji}\}$ to a (pro-)object $F$ of $\hat{\mc{U}}$. 

Since $\phi_i^*(\omega_j)\in p\cdot H^0(U_i,\Omega^1(D))$ for every $j$, which is evident if $j\neq i$, and follows from $\phi_i^*(I_{a_i})=I_{a_i}^p$ if $j=i$, we conclude the existence of the induced connection on $F_{i}$. In order to show that $\{F_{i}\}_{i=0,\dots,s}$ glues, it suffices to show
\begin{equation}\label{equation-gluing-to-show}
g_{ji}(p^{[r]}\phi_i^*(\alpha))\in F_{j}(U_i\cap U_j)
\end{equation}
for all $\alpha\in \Gamma^{\otimes r}$. In fact, $g_{ji}(p^{[r]}\phi_i^*(\alpha))=p^{[r]}\cdot g_{ji}(\phi_i^*(1))\cdot \alpha$, because $g_{ji}$ is compatible with the right multiplication on $E$ by $W(k)\fps{\Gamma}$. Therefore $p^{[\ell]}\cdot p^{[r]}\subset p^{[\ell+r]}$ reduces \eqref{equation-gluing-to-show} to the case $\alpha=1$. Since the image of $\nabla_{\partial_t}^m(1)$ via the projection to $\OO_{U_i\cap U_j}\otimes_{W(k)} \Gamma^{\otimes r}$ vanishes if $m<r$, we are done.
  
We denote by $\epsilon$ the inclusion $(F,\nabla)\xr{} F^*(E,\nabla)$. Once again, for $x\in \{(a_1,\xi),y\}$, let us identify 
$$
\omega_{x}(F^*(E,\nabla)) = W(k)\hat{\otimes}_{\sigma,W(k)} W(k)\fps{\Gamma} 
$$ 
via $\epsilon_1$, so that 
\begin{equation}\label{equation-use-of-F}
\omega_{x}(\epsilon)(\omega_{x}(F,\nabla))=W(k)\hat{\otimes}_{\sigma,W(k)} \prod_{r\geq 0} p^{[r]} \cdot \Gamma^{\otimes r}.
\end{equation}
It follows easily from the functoriality of $\eta_x$ and its simple shape for $\mathbf{1}$, that 
\begin{align*}
\eta_x(1)&\in 1+W(k)\hat{\otimes}_{\sigma,W(k)} W(k)\fps{\Gamma}\cdot \Gamma,\\
\eta_x(1\otimes \alpha)&=\eta_x(1)\cdot (1\otimes \alpha),\quad \text{for all $\alpha\in W(k)\fps{\Gamma}$.}
\end{align*}
We claim that, for all $m\geq 0$, 
\begin{equation}\label{equation-eta-maps}
\eta_x(1\otimes \Gamma^{\otimes m}) \subset W(k)\hat{\otimes}_{\sigma,W(k)}  \prod_{r\geq 0} p^{[r]} \cdot \Gamma^{\otimes r+m}.
\end{equation}
Suppose first that $x=y$, and fix a coordinate $t$ on $U_1$. Then $\eta_y$ is the $W(k)\hat{\otimes}_{\sigma,W(k)}$-base change of 
$$
\alpha\mapsto \sum_{m=0}^{\infty}\frac{(t(y)-t(\phi(y)))^m}{m!} \cdot \nabla_{\partial_t}^m(\alpha)_{\mid t=\phi(y)}, 
$$
and we can argue as above. If $x=(a_1,\xi)$ then $\eta_{(a_1,\xi)}$ is the base change of 
$$
\alpha\mapsto \exp\left( \log(\lambda_{\xi,\phi}) \cdot {\rm Res}_{a_1}(E,\nabla)\right)(\alpha),
$$ 
for a certain $\lambda_{\xi,\phi}\in W(k)^*$ (Section \ref{section-tangential-base-points}), and the residue ${\rm Res}_{a_1}(E,\nabla)\in {\rm End}(\omega_{a_1}(E))$. Since $\log(\lambda_{\xi,\phi})\in p\cdot W(k)$ and$$
{\rm Res}_{a_1}(E,\nabla)(\alpha)=\omega_1^{\vee}\cdot \alpha,
$$
we can even show \eqref{equation-eta-maps} with substituting $\frac{p^r}{r!}$ for $p^{[r]}$.

In view of the bijection
\begin{equation}\label{equation-universal-property}
\Hom_{\hat{\mc{U}}}(E,F)\xr{\cong} \omega_{a_1}(F), \quad f\mapsto \omega_{a_1}(f)(1), 
\end{equation} 
for any $F\in {\rm ob}(\hat{\mc{U}})$, and \eqref{equation-eta-maps}, there is a unique $\theta':(E,\nabla)\xr{} (F,\nabla)$ such that $\omega_{a_1}(\theta')(1)=\eta_{(a_1,\xi)}(1)$. We get a commutative diagram
$$
\xymatrix{
\omega_{a_1}(F^*(E)) \ar[rr]^{\gamma_{y,(a_1,\xi),F^*(E)}} 
&
&
\omega_{y}(F^*(E))
\\
\omega_{a_1}(F) \ar[rr]^{\gamma_{y,(a_1,\xi),F}} \ar[u]^{\omega_{a_1}(\epsilon)}
&
&
\omega_{y}(F) \ar[u]_{\omega_{y}(\epsilon)}
\\
 \omega_{a_1}(E)  \ar[rr]^{\gamma_{y,(a_1,\xi),E}} \ar[u]^{\omega_{a_1}(\theta')}
&
&
\omega_{y}(E), \ar[u]_{\omega_{y}(\theta')} 
}
$$     
which, together with the fact that $\gamma_{y,(a_1,\xi)}$ is an isomorphism of $F$-fibre functors, implies 
$$ 
\eta_{y}(1\otimes \gamma_{y,(a_1,\xi),E}(1))= \omega_y(\epsilon)(\omega_y(\theta')(\gamma_{y,(a_1,\xi),E}(1))).
$$ 
In view of \eqref{equation-use-of-F} and \eqref{equation-eta-maps}, we conclude \eqref{equation-to-prove-in-main-thm}.
\end{proof}
\end{thm}

\subsection{}
Suppose that $a_0=\infty,a_1=0,a_2=1$; we set $A=\omega_1^{\vee}$ and $B=\omega_2^{\vee}$. Furusho's $p$-adic multiple zeta values can be computed as follows. Let $\xi_0$ be the tangent vector at $0$ pointing to $1$, and let $\xi_1$ be the tangent vector at $1$ pointing to $0$. We denote by $x_0=(\omega_0,\eta_{\xi_0})$ and $x_1=(\omega_1,\eta_{\xi_1})$ the corresponding tangential $F$-fibre functors. As noted in Section \ref{section-tangential-Deligne}, they do not depend on a choice of $\log(p)$. 

		Furusho \cite[Theorem~3.3/Proposition~3.7]{F} constructs solutions $G^a_0(y)$ and $G^a_{1}(y)$ of the $p$-adic Knizhnik-Zamolodchikov equation with values in Coleman functions. In view of Proposition \ref{proposition-comparison}, it is easy to see that, for every $y\in W(k)$ with $\bar{y}\not\in \{0,1\}$, and every finite extension $k\supset \F_p$,
		\begin{equation}\label{equation-path-in-two-steps}
		G^a_0(y)=\gamma_{y,x_0}(1), \qquad 	G^a_1(y)=\gamma_{y,x_1}(1).
		\end{equation}
		Furusho defines the $p$-adic Drinfel'd associator $\Phi^p_{KZ}(A,B)$ by requiring 
		$$
		G^a_0(y)=G^a_1(y)\cdot \Phi^p_{KZ}(A,B).
		$$
                Since right multiplication by $W(k)\fps{\Gamma}$ induces endomorphisms of the pro-object $\hat{E}$, we have
                \begin{equation}\label{equation-right-mult} 
		\gamma_{z,x}(e)=\gamma_{z,x}(1)\cdot e
		\end{equation}
		for every two $F$-fibre functors $x,z$. Therefore
		\begin{equation}\label{equation-product-formulas}
		\gamma_{x_1,y}(1)\cdot \gamma_{y,x_1}(1)=\gamma_{x_1,x_1}(1)=1,\qquad 
		\gamma_{x_1,y}(1)\cdot \gamma_{y,x_0}(1)=\gamma_{x_1,x_0}(1),
		\end{equation}
hence $\Phi^p_{KZ}(A,B)=\gamma_{x_1,x_0}(1)$. Now it follows from Furusho's definition of $p$-adic multiple zeta values as coefficients of the Drinfel'd associator that 
\begin{equation}\label{equation-def-mzv}
(-1)^m \zeta_p(k_1,\dots,k_m)=c(k_1,\dots,k_m) 
\end{equation}
where $c(k_1,\dots,k_m)$ is the coefficient of $A^{k_m-1}BA^{k_{m-1}-1}B\dots A^{k_1-1}B$ in $\gamma_{x_1,x_0}(1)$.

\begin{corollary}\label{corollary-multiple-zeta-values}
For all $k_1,\dots,k_m$, we have $\zeta_p(k_1,\dots,k_m)\in p^{\left[\sum_{i=1}^m k_i\right]}$.
\begin{proof}
Follows immediately from \eqref{equation-def-mzv}, \eqref{equation-product-formulas}, and Theorem \ref{thm-main}. 
\end{proof}
\end{corollary}

\bibliography{pzeta}
\end{document}